\documentclass[11pt]{amsart}
\usepackage{amsmath,amssymb}
\numberwithin{equation}{subsection}
\theoremstyle{plain}
  \newtheorem{theorem}{Theorem}[section]
  
  \newtheorem{proposition}[theorem]{Proposition}
  \newtheorem{lemma}[theorem]{Lemma}
  \newtheorem{corollary}[theorem]{Corollary}

\theoremstyle{definition}
  \newtheorem{definition}[theorem]{Definition}
  \newtheorem{example}[theorem]{Example}
\theoremstyle{remark}
 \newtheorem{remark}[theorem]{Remark}

\newcommand{\brac}[2]{{\langle #1\mid#2\rangle}}

\newcommand{\Ch}[1]{ {\mathrm{Tr}} ( {#1}; q^\partial)}

\newcommand{\Worb}[1]{{}_{\W}\!{#1}}
\newcommand\qfac[1]{[{#1}]_q!}

\newcommand\projmapop[1]{\projmap({#1})}

\def\r{r}

\def\C{{\mathbb F}}
\def\Q{{\mathbb Q}}
\def\Z{{\mathbb Z}}
\def\F{{\mathbb F}}

\def\Qgeq{\F\setminus\Q_{\leq0}}

\def\projmap{\varrho}
\def\diag{\theta}
\def\rat{{}}

\def\af{{\rm aff}}
\def\peri{{\gamma}}

\def\ratH{{\H_\kappa}}
\def\Preg{P_\circ}
\def\affH{{\H^\af}}

\def\O{{\mathcal O}}

\def\tab{\mathop{\mathbf{tab}}\nolimits}
\def\tabRC{{\mathbf{st}}}
\def\Tab{{\mathop{\mathbf{Tab}}\nolimits}}

\def\PP{{\mathbf{PP}}}
\def\PPR{{\PP}^{\mathrm R}}
\def\PPC{{\PP}^{\mathrm C}}
\def\trigPP{\widetilde{\PP}}
\def\trigPPR{{\trigPP}{}^{\mathrm R}}
\def\trigPPC{{\trigPP}{}^{\mathrm C}}

\def\PPp{\PP_\pri}
\def\PPRp{{\PP}_\pri^{\mathrm R}}
\def\PPCp{{\PP}_\pri^{\mathrm C}}
\def\trigPPp{{\trigPP}{}_\pri}
\def\trigPPRp{{\trigPP}{}_\pri^{\mathrm R}}
\def\trigPPCp{{\trigPP}{}_\pri^{\mathrm C}}
\def\GF{{\Psi}}

\def\GFR{\GF(\PPR(\diag);q)}
\def\GFRperi{\GF(\PPRp(\diag);q)}
\def\GFC{\GF(\PPC(\diag);q)}
\def\GFCperi{\GF(\PPCp(\diag);q)}

\def\mV{{\mathcal{V}}}
\def\pri{\peri}
\def\TabRC{{\mathbf{St}}_\pri}
\def\ratL{{\mathcal L}_\kappa^\rat}
\def\Vlm{{\trigVdiag}}
\def\trigVdiag{\trigV(\adiag)}
\def\affV{{\mV^\af}}

\def\ratV{{\ratL}}
\def\ratVlm{{\ratV(\lm)}}
\def\polyVlm{{\widetilde{\mV}_\pri(\alm)_{\posi}}}
\def\trigV{\widetilde{\mV}_\pri}

\def\TabRCplm{{\TabRC(\alm)_{\posi}}}
\def\st{\tabRC}

\def\H{{\mathcal H}}
\def\trigH{{\widetilde\H_\kappa}}
\def\trigHm{{\widetilde\H_{-\kappa}}}

\def\trigHmod{{\trigH{\hbox{-}}{\mathrm{mod}} }}

\def\ratst{{\Delta}_\kappa^\rat}
\def\ratHmod{{ \ratH{\hbox{-}}{\mathrm{mod}} }}
\def\ratH{{\H_\kappa}}

\def\affH{{\H^\af}}

\def\affVlsm{{\affV}(\diag)}
\def\polyH{{\H^\circ_{\kappa}}}
\def\posi{+}
\def\adiag{\hat{\diag}}
\def\restab{\tabRC}
\def\WRC{{\mathcal Z}}
\def\mkappa{(-m,\kappa-m)}

\def\Ind{\mathbf{ind}}
\def\nil{\mathbf{nil}}

\def\Irr{{\hbox{\rm{Irr}}}}
\def\dim{{\hbox{\rm{dim}}}}
\def\id{{\hbox{\rm{id}}}\,}
\def\mod{{\rm mod\,}}

\def\al{{\alpha}}
\def\alch{\alpha^\vee}
\def\e{{\epsilon}}

\def\ech{\epsilon^\vee}
\def\lm{{\lambda}}

\def\lsm{{\lambda/\mu}}
\def\alm{\wh{\lm}}

\def\m{{m}}
\def\gl{{\mathfrak{gl}}}
\def\h{{\mathfrak h}}

\def\bx{\framebox(10,10)}

\def\bra{{\langle}}
\def\ket{{\rangle}}
\def\qed{\hfill$\square$}
\def\udl{\underline}

\def\wh{\widehat}
\def\+{\mathop{\oplus}}
\def\*{\mathop{\otimes}}
\def\isomto{\stackrel{\sim}{\to}}
\def\O{{\mathcal O}}
\def\Oss{{\mathcal O}^{\h}}

\def\XX{{\F[\underline{x}^{\pm1}]}}
\def\YY{{\F[\underline{y}]}}

\def\xch{y}
\def\Xch{\F[\underline{y}]}
\def\X{\F[\underline{x}]}
\def\XL{\F[\underline{x}^{\pm1}]}
\def\UU{S(\h)}

\def\xch{y}
\def\sp{S}

\def\emb{\iota}
\def\embinv{\jmath}

\def\part{{\Lambda}^+}
\def\affWreg{{\affW}_\circ}

\def\part{{\Lambda}^+}

\def\aa{a}
\def\W{{\mathcal{W}}}
\def\Weyl{{\W}}
\def\affW{ { \wh\W } }

\def\Zlm{Z^{\alsm}}
\def\Zlmp{{{Z}^{\wh\lm}_\circ}}

\def\dom{{\mathfrak P}^+}

\def\rbtab{{\mathfrak{t}_\diag }}

\def\one{{\bf 1}}
\def\con{{\mathfrak c}}

\def\weight{\dot\con}
\def\t{{t}}
\def\dset{{\mathbf{Y}}}

\def\dsetint{{}_\Z\!\dset}
\def\unip{{\bf e}}

\def\vert{\delta_n}
\def\avert{\wh\delta_n}
\def\fdiag{{{\theta}}}

\def\der{\partial}

\def\cas{d}

\def\p{{\mathfrak{p}}}
\def\vec{{\mathrm{v}}}
\begin{document}
\title[Cherednik Algebras]{
Cylindrical Combinatorics and Representations of
Cherednik Algebras of type $A$}
\author{Takeshi Suzuki}
\address{Research Institute for Mathematical Sciences\\
Kyoto University, Kyoto, 606-8502, Japan}
\email{takeshi@kurims.kyoto-u.ac.jp}
\begin{abstract}
We investigate the representation theory of
the rational and trigonometric Cherednik algebra of type $GL_n$
by means of combinatorics on periodic (or cylindrical)
skew diagrams.

We introduce and study standard tableaux and plane partitions
on periodic diagrams, and
in particular, compute some generating functions 
concerning plane partitions, where 
Kostka polynomials and their level restricted generalization
appear.

On representation side, we study representations of
Cherednik algebras which admit weight decomposition 
with respect to a certain commutative subalgebra.
All the irreducible representations of this class
are constructed combinatorially 
using standard tableaux on periodic diagrams,
and this realization as "tableaux representations"
provides a new combinatorial approach to the investigation of
these representations.

As consequences, we describe the decomposition of a tableaux representation
as a representation of the degenerate affine Hecke algebra,
which is a subalgebra of the Cherednik algebra,
and also describe
the spectral decomposition of the spherical subspace
(the invariant subspace under the action of the Weyl group)
of a tableaux representation
with respect to the center of the degenerate affine Hecke algebra,
In particular, the computation of the character
of the spherical subspace is reduced to 
the computation of the generating function
for the set of column strict plane partitions,
and we obtain an expression of the characters 
in terms of Kostka polynomials as announced in \cite{Su;conformal}.
\end{abstract}
\maketitle

%\bigskip
\noindent
{\bf CONTENTS}

\smallskip\noindent
\S1 Plane partitions and standard tableaux on periodic diagrams

\smallskip\noindent
\S2 Generating functions

\smallskip\noindent
\S3 Tableau representations of the trigonometric Cherednik algebra

\smallskip\noindent
\S4 Application to the rational Cherednik algebra

\smallskip\noindent
\S5 Characters 

\smallskip\noindent
Appendix A. Classification of irreducible modules with weight decomposition
%%%%%%%%%%%%%%%%%%%%%%%%%%%%%%%%%%%%%%%%%%%%%%%%%%%%%%%%
\section{Plane partitions and standard tableaux 
on periodic diagrams}
%%%%%%%%%%%%%%%%%%%%%%%%%%%%%%%%%%%%%%%%%%%%%%%%%%%%%%%%%%%%
We study cylindrical combinatorics, that is, 
combinatorics on 
{\it periodic $($or cylindrical$)$ skew diagrams},
which are introduced by Gessel and Krattenthaler \cite{GK}
as a cylindrical analogue of skew Young diagrams,
and have appeared in the representation theory of 
the double affine Hecke algebras \cite{Ch;fourier,SV}.
%
%%%%%%%%%%%%%%%%%%%%%%%%%%%%%%%%%%%%%%%
\subsection{Root system and Weyl group}
%%%%%%%%%%%%%%%%%%%%%%%%%%%%%%%%%%%%%%
Let  $\F$ denote a field of characteristic $0$, which
 includes the field $\Q$ of rational numbers
and the ring $\Z$ of integers.

Throughout this article, we use the following notation:

$$[a,b]=
\begin{cases}
\{a,a+1,\dots,b\}&\hbox{ for }a,b\in\F
\hbox{ with }b-a\in\Z_{\geq0},\\
\emptyset&\hbox{ otherwise}.
\end{cases}
$$
%%%%%%%%%%%%%%%%%%%%%%%%%%%%%%%%%%%%%%%%%%%%%%%

Fix $n\in\Z_{\geq 2}.$
Let  $\h$ be an $n$-dimensional
vector space over $\F$
with the basis $\{\ech_1,\ech_2,\dots,\ech_n\}$:
$$\h=\oplus_{i\in[1,n]}\F\ech_i.$$
Introduce the non-degenerate  symmetric bilinear form $(\ |\  )$
on $\h$ by
$(\ech_i|\ech_j)=\delta_{ij}.$
Let $\h^*=\oplus_{i=1}^n\F\e_i$
be the dual space of $\h$,
where $\e_i$ 
are  the dual vectors of $\ech_i$.
The natural pairing is denoted by
$\bra\, |\,\ket : \h^*\times \h\to\F$.

Put $\al_{ij}=\e_i-\e_j$,
$\alch_{ij}=\ech_i-\ech_j\ (1\leq i\neq j\leq n)$
and $\al_i=\e_i-\e_{i+1},\
\alch_i=\ech_i-\ech_{i+1}\ (i\in[1,n-1])$.
Then
\begin{align*}
R& =\left\{ \al_{ij}\mid i,j\in[1,n],\ i\ne j \right\},\ \
R^+=\left\{ \al_{ij}\mid i,j\in[1,n],\ i<j\right\}
\end{align*}
give the system of roots and
positive roots
of type $A_{n-1}$ respectively.

%%%%%%%%%%%%%%%%%%%%%%%%%%%%%%%%%%%%%%
Let $\W$ %and $\afW$
denote the Weyl group associated with the
root system $R$. %and $\wh R$ respectively.
The group $\W$  acts on $\h$ and $\h^*$, and
it is isomorphic to the symmetric group
of degree $n$.

Denote by $s_\al$
the reflection in $\W$ corresponding to
$\al\in  R$.
We write
 $s_{ij}=s_{\al_{ij}}$ $(i,j\in[1,n])$ and $s_i=s_{\al_i}$
$(i\in[1,n-1])$.
We have $\W=\bra s_1,s_2,\dots,s_{n-1}\ket$.
%%%%%%%%%%%%%
%

Put
$$P=\oplus_{i\in[1,n]} \Z\e_i,
$$
which is a subset of $\h^*$ and is preserved by the action of $\W$.
Define the {\it extended affine Weyl group} $\affW$ of $\gl_n$
as the semidirect product $P\rtimes \W$.

For $\eta\in P$, let  $\t_\eta$ denotes the corresponding
 element in $\affW$.
Put $s_0=\t_{\e_1-\e_n}s_{1n}$ and
$\pi=\t_{\e_1}s_{1}s_2\dots s_{n-1}$.
Then $\{ s_0,s_1,\dots,s_{n-1},\pi^{\pm1}\}$ gives a set of
generators of $\affW$,
and the subgroup generated by $s_0,s_1,\dots,s_{n-1}$
is the affine Weyl group associated
with the root system of type $A_n^{(1)}$.

Define an action of the extended affine Weyl group $\affW$ on the set $\Z$
of integers by
\begin{align}
\label{eq;WactZ1}
s_i(j)&=
\begin{cases}
j+1&\quad\hbox{for }j\equiv i\ \mod n, \\
j-1&\quad\hbox{for }j\equiv i+1\ \mod n,\\
j&\quad\hbox{for }j\not\equiv i,i+1\ \mod n,
\end{cases}
\ \ (i\in[1,n-1])\\
\t_{\e_i}(j)&=
\begin{cases}
j+n\quad&\hbox{for }j\equiv i\ \mod n,\\
j\quad&\hbox{for }j\not\equiv i\ \mod n.
\end{cases}
\ \ (i\in[1,n])
\end{align}
Note in particular that $\pi(j)=j+1$ for all $j\in\Z$.
%%%%%%%%%%%%%%%%%%%%%%%%%%%%%%%%%%%%%%%%%%%%%%%%%%%%%
\subsection{Periodic skew diagrams}
%%%%%%%%%%%%%%%%%%%%%%%%%%%%%%%%%%%%%%%%%%%%%%%%%%%

We need a slightly generalized definition of skew diagrams.
%%%%%%%%%%%%%%%%%%%%%%%%%%%%%%%%%%%%%%%%%%%%%%%%%%%%%%%%%%%%
\begin{definition}\label{def;skew_diagram}
A subset $\diag$
of $\Z\times\F$
is called a {\it skew diagram}
%of order $n$
if the following conditions are satisfied:

\smallskip\noindent
(D1) The set $\diag$ consists of finitely many elements.

\smallskip\noindent
(D2) For any $i\in\Z$, there exist
$\lm_i,\mu_i\in\F$ such that
 $\{(a,b)\in\fdiag\mid a=i\}
=[\mu_i+1,\lm_i]$.

\smallskip\noindent
(D3) (Skew property)
If $(a,b),(a',b')\in \fdiag$
with  $a'-a\in\Z_{\geq1}$
and $(b'-a')-(b-a)\in\Z_{\geq 0}$ then
$(a,b+1),(a',b'-1)\in\fdiag$.
\end{definition}
%%%%%%%%%%%%%%%%%%%%%%%%%%%%%%%%%%%%%%%%%%%%%%%%%%%%%%

In the sequel, we regard $\Z\times \F$ as a $\Z$-module.

Let $\peri\in\Z\times\F$. Denote by $\Z\peri$ the subgroup of $\Z\times\F$
generated by $\peri$.
%%%%%%%%%%%%%%%%%%%%%%%%%%%%%%%%%%%%%%%%%%%%%%%%%%%%%%%%%%%%%%%%%%%
\begin{definition}\label{def;periodic_skew}
A subset $\Theta$  of $\Z\times\F$
is called a {\it  periodic skew diagram}
(or a {\it cylindrical skew diagram})
of period $\peri$
if the following conditions are satisfied:

\smallskip\noindent
(D'1)  The group $\Z\peri$ acts on $\Theta$ by
 parallel translation, i.e., $\Theta+\peri=\Theta$,
and
a fundamental domain of this action on $\Theta$
consists of finitely many elements.

\smallskip\noindent
(D'2) For any $i\in\Z$, there exist
$\lm_i,\mu_i\in\F$ such that
 $\{(a,b)\in\Theta\mid a=i\}
=[\mu_i+1,\lm_i]$.

\smallskip\noindent
(D'3) (Skew property)
If $(a,b),(a',b')\in \Theta$
with  $a'-a\in\Z_{\geq1}$
and $(b'-a')-(b-a)\in\Z_{\geq 0}$ then
$(a,b+1),(a',b'-1)\in\Theta$.
\end{definition}
%%%%%%%%%%%%%%%%%%%%%%%%%
Let $\dset^n$ denote the set of skew diagram
consisting of $n$-elements, and
let $\wh\dset_\peri^n$ denote the set of
periodic diagrams of period $\peri$ consisting of $n$ numbers of
$\Z\peri$-orbits.

%%%

For a skew diagram  $\diag$, 
define a subset of $\Z\times\F$ by
$$\adiag_\peri=\diag+\Z\peri.$$
We often write $\adiag=\adiag_\peri$ when $\peri$ is fixed.
For $\peri=(\pm m,l)\in\Z\times\F$ with $m\in\Z_{\geq1}$,
set
$$\dset^n_\peri=\{\diag\in\dset^n\mid
\diag\subset [1,m]\times\F\hbox{ and }
\adiag_\peri\in\wh\dset_\peri^n \}.$$
Note that $\wh\dset^n_\peri=\wh\dset^n_{-\peri}$ and
$\dset^n_\peri=\dset^n_{-\peri}$.

For a periodic skew diagram $\Theta$ of period $\peri=(\pm m,l)$,
the subset $\Theta\cap ([1,m]\times\F)$ is a skew diagram and
it gives a fundamental domain of $\Z\peri$ on $\Theta$.
Hence any periodic diagram of period $\peri$
is of the form $\adiag$ for some $\diag\in\dset^n_\peri$,
and the map $\diag\mapsto\adiag$,
$\dset_\peri^n\to\wh\dset_\peri^n$
is bijective.

By the condition (D2), any skew diagram
$\diag\subset [1,m]\times\F$ is expressed %up to parallel shift
as $\lm/\mu$ for some $\lm,\mu\in\F^m$, % and some $m\in\Z_{\geq1}$,
where
$$\lm/\mu=\{
(a,b)\in\Z\times\F\mid a\in[1,m],\ b\in[\lm_a+1,\mu_a]\}.$$

\def\domset{{\mathfrak D}}
Moreover, it follows from the condition (D3) that
$\lm$ and $\mu$ can be chosen from the set of dominant elements
$$\domset_m=\{\nu=(\nu_1,\dots,\nu_m)\in\F^m\mid
(\nu_i-i)-(\nu_j-j)\notin\Z_{\leq0}\ \hbox{ for any }i<j\}.
$$

For $\kappa\in\F$, define
$\domset_{m,\kappa}$ as the subset of $\domset_m$ consisting of those 
elements
$\nu$ satisfying the following conditions:
\begin{align}\label{eq;domsetmk1}
&p\kappa+(\nu_i-i)-(\nu_j-j)\notin\Z_{\leq0}\ \hbox{ for any }i<j
\hbox{ and }p\in\Z_{\geq0},\\
\label{eq;domsetmk2}
&p\kappa-(\nu_i-i)+(\nu_j-j)\notin\Z_{\leq0}\ \hbox{ for any }
i<j \hbox{ and }p\in\Z_{>0}.
\end{align}
Note that
$\domset_{m,\kappa}=\emptyset$ unless $\kappa\in\F\setminus\Q_{\leq0}$.
For $\nu\in\Z^m$,
we write $\nu\models n$ if $\nu$ is a composition of $n$; namely,
$\nu_1\geq \nu_2\geq\dots\geq\nu_m$ and 
$\sum_{i\in[1,m]}\nu_i=n$.

It is easy to show the following:

%%%%%%%%%%%%%%%%%%%%%%%%%%%%%%%%%%%%%%%%%%%%%%%%%%%%%%%%%%%
\begin{lemma}\label{lem;dominant}
 Let $m\in\Z_{\geq1}$ and 
$\kappa\in\F\setminus\Q_{\leq0}$.
If $\lm,\mu\in\domset_{m,\kappa}$ and $\lm-\mu\models n$,
 then $\lsm\in\dset^n_{(-m,\kappa-m)}$.
Conversely, for any $\diag\in\dset_{(-m,\kappa-m)}^n$,
there exist
 $\lm,\mu\in\domset_{m,\kappa}$ such that 
$\diag=\lsm$.
\end{lemma}

Define 
$\dset_\peri^{*n}$ as the subset of $\dset_\peri^n$ 
consisting of diagrams without empty rows:
$$\dset_\peri^{*n}=
\left\{\diag\in\dset_\peri^n\;\left|\;
\forall a\in\Z,\  \exists b\in\F \hbox{ such that }
(a,b)\in\adiag\right.\right\}.$$
%%%%%%%%%%%%%%%%%%%%%%%%%%%%%%%%%%%%%%%%%%%%%
\begin{corollary}\label{cor;dominant}
Let $m\in[1,n]$ and
$\kappa\in\F\setminus\Q_{\leq0}$.
Put $\peri=(-m,\kappa-m)$.

\smallskip\noindent
$\mathrm{(i)}$
Let $\lm,\mu\in\F^m$ such that $\lm-\mu\models n$
and $\lm_i-\mu_i\geq1$ for all $i\in[1,m]$.
Then $\lm/\mu$ is in $\dset_\peri^{*n}$ if and only if
$\lm,\mu\in\domset_{m,\kappa}$.

\smallskip\noindent
$\mathrm{(ii)}$
For any $\diag\in\dset_\peri^{*n}$,
there exist unique $\lm,\mu\in\domset_{m,\kappa}$ such that 
$\diag=\lsm$.
\end{corollary}
%%%%%%%%%%%%%%%%%%%%%%%%%%%%%%%%%%%%%%%%%%
%
%%%%%%%%%%%%%%%%%%%%%%%%%%%%%%%%%%%%%%%%%%%%%%%%%%%%%%%%
\subsection{Plane partitions on periodic diagrams}
%%%%%%%%%%%%%%%%%%%%%%%%%%%%%%%%%%%%%%%%%%%%%%%%%%%%%%
%

Let $\diag$ be a skew diagram.
\begin{definition}
(i)  A map $\p:\diag\to\Z$
is called {\it a plane partition} on $\diag$
if it is weakly row-column increasing; namely, if it
satisfies the following two conditions (P2) and (P3):

\smallskip\noindent
$\mathrm{(P2)}\ \p(a,b)\leq \p(a,b+1)$ for any
$(a,b)\in\diag$ and  $(a,b+1)\in\diag$.
%with $k\in\Z_{\leq-1}$.

\smallskip\noindent
$\mathrm{(P3)}\ \p(a,b)\leq \p(a+k+1,b+k)$ for any
 $(a,b)\in\diag$ and $(a+k+1,b+k)\in\diag$ %from $\diag$
with $k\in\Z_{\geq0}$.

\smallskip\noindent
(ii) A plane partition on $\diag$ is said to be {\it row strict}
(resp., {\it column strict})
if the strict inequality always holds in (P2)
 (reps., (P3)).
\end{definition}
Let $\trigPP(\diag)$ denote the set of
the plane partitions on $\diag$, and
let $\trigPPR(\diag)$ and $\trigPPC(\diag)$
 denote the set of the row strict and column strict
plane partitions on $\diag$ respectively.
Define
\begin{align}
&\PP(\diag)%_{\geq0}
=\{\p\in\trigPP(\diag)\mid \p(u)\geq 0\ \  \forall u\in\diag\},
\end{align}
and define $\PPR(\diag)$ and  $\PPC(\diag)$
similarly.

\begin{remark}
In literature $($e.g. \cite{St2}$)$,
a plane partition on $\diag\subset\Z\times\Z$
is defined as a map $\p:\diag\to\Z_{\geq1}$
%(or $\p:\diag\to\Z_{\geq1}$)
which is weakly decreasing in both row and column directions.
In \cite{St2}, an element of $\PP(\diag)$
is called a weak reverse plane partition.
\end{remark}
Let $\peri\in\Z\times\F$.
%=\mkappa$ with $m\in\Z_{\geq1}$
%and $\kappa\in\F\setminus\{0\}$.
%$\Q_{\leq0}$.
Let $\diag\in\dset_\peri^n$.
%%%%%%%%%%%%%%%%%%%%%%%%%%%%%%%%%%%%%%%%%%%%%%%%%%%%
\begin{definition}\label{def;PP}
(i)  A map $\p:\adiag\to\Z$
is called {a plane partition} on $\adiag$
%(or,  a $\peri$-{\it periodic plane partition} on $\diag$)
if it satisfies the following conditions:

\smallskip\noindent
(P1) $\p(u+\peri)=\p(u)-1$ for all $u\in\adiag$.
%%perisgn

\smallskip\noindent
(P2) $\p(a,b)\leq \p(a,b+1)$ for any
$(a,b)\in\adiag$ and $(a,b+1)\in\adiag$. %from $\adiag$.
%with $k\in\Z_{\leq-1}$.

\smallskip\noindent
(P3) $\p(a,b)\leq \p(a+k+1,b+k)$ for any
 $(a,b)\in\adiag$ and $(a+k+1,b+k)\in\adiag$ %from $\adiag$
with $k\in\Z_{\geq0}$.

\smallskip\noindent
(ii) A plane partition on $\adiag$ is said to be {\it row strict}
(resp., {\it column strict})
if the strict inequality always holds in the condition
(P2) (resp., (P3)).
\end{definition}
Let $\trigPPp(\diag)$ denote the set of
the plane partitions on $\adiag$, and
let $\trigPPRp(\diag)$ (resp., $\trigPPCp(\diag))$
 denote the set of the row strict (resp., column strict)
plane partitions on $\adiag$.
Define
$$\PPp(\diag)=
\{\p\in\trigPPp(\diag)\mid \p(u)\geq 0\ \  \forall u\in\diag\},$$
and define $\PPRp(\diag)$ and $\PPCp(\diag)$
similarly.
\begin{example}
Let $n=7$, $\peri=(-2,3)$, and let
%$$\diag=\{(1,2),(1,3),(1,4),(1,5),(2,1),(2,2),(2,3)\};$$
%namely,
 $\diag=\lm/\mu$ with $\lm=(5,3)$ and $\mu=(1,0)$.
Then $\diag\in\dset_\peri^n$.
The following figure represents
the associated periodic diagram $\adiag$
and a column strict plane partition on $\wh\diag$.
%%%%%%%%%%%%%%%%%%%%%%%%%%

\smallskip
\begin{center}
\setlength{\unitlength}{.06cm}
\begin{picture}(150,90)(-40,-40)
\multiput(65,40)(5,0){5}{\rule{.8pt}{.8pt}}
\put(40,20){\bx{$-\!1$}}
\put(50,20){\bx{$0$}}
\put(60,20){\bx{$0$}}
\put(70,20){\bx{$2$}}
\put(30,10){\bx{$0$}}
\put(40,10){\bx{$1$}}
\put(50,10){\bx{$1$}}
\put(10,0){\bx{$0$}}
\put(20,0){\bx{$1$}}
\put(30,0){\bx{$1$}}
\put(40,0){\bx{$3$}}
\put(0,-10){\bx{$1$}}
\put(10,-10){\bx{$2$}}
\put(20,-10){\bx{$2$}}
\put(85,-3){$\;\diag$}
\put(-20,-20){\bx{$1$}}
\put(-10,-20){\bx{$2$}}
\put(0,-20){\bx{$2$}}
\put(10,-20){\bx{$4$}}
\put(-30,-30){\bx{$2$}}
\put(-20,-30){\bx{$3$}}
\put(-10,-30){\bx{$3$}}
\multiput(-35,-40)(5,0){5}{\rule{.8pt}{.8pt}}
\put(80,-3)
%{\mbox{$\}$}}
{\mbox{$\biggl. \biggr\}$}}
\linethickness{1.5pt}
\put(10,10){\line(1,0){40}}
\put(10,10){\line(0,-1){10}}
\put(0,0){\line(1,0){10}}
\put(0,0){\line(0,-1){10}}
\put(0,-10){\line(1,0){30}}
\put(30,0){\line(1,0){20}}
\put(30,0){\line(0,-1){10}}
\put(50,0){\line(0,1){10}}
%extra border
\end{picture}

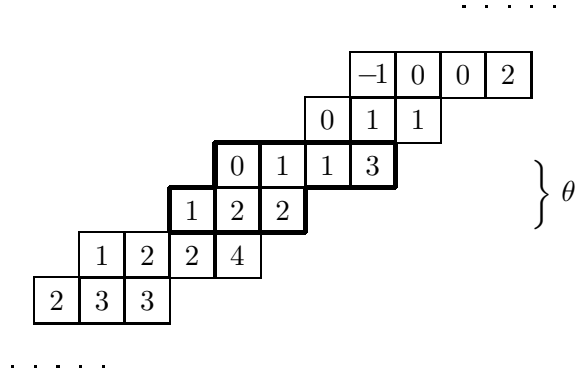
{\footnotesize {Figure 1.}}
\end{center}
\end{example}
\def\Map{{\mathrm{Map}}}

%%%%%%%%%%%%%%%%%%%%%%%%
A plane partition $\p\in\trigPP(\diag)$ on $\diag$
can be uniquely extended to a function $\hat\p:\adiag\to \Z$
by setting
$\hat\p(u+k\peri)=\p(u)-k$ for $u\in\diag$ and $k\in\Z$.
The correspondence $\p\mapsto\hat\p$ gives an embedding
from $\trigPP(\diag)$ into the set $\Map(\adiag,\Z)$
of maps $\adiag\to\Z$.
% satisfying (P1).

In the sequel, 
we often identify $\trigPP(\diag)$
with a subset of $\Map(\adiag,\Z)$.
Under this identification, we have
$$\trigPPp(\diag) %\resd
\subseteq \trigPP(\diag),\
\trigPPRp(\diag) %\resd
\subseteq \trigPPR(\diag),\
\trigPPCp(\diag)%\resd
\subseteq \trigPPC(\diag)$$
$$\PPp(\diag)%\resd
\subseteq\PP(\diag),\ 
\PPRp(\diag) %\resd
\subseteq \PPR(\diag),\ 
\PPCp(\diag)%\resd
\subseteq \PPC(\diag).$$
%Similarly to Lemma \ref{lem;equi_TabRC}, we have

The following statement can be shown easily using
the skew property of $\adiag$:
%%%%%%%%%%%%%%%%%%%%%%%%%%
\begin{lemma}
  \label{lem;equi_PP}
 A plane partition
 $\p\in\trigPP(\diag)$ is in $\trigPP_\pri(\diag)$
if and only if
$$\p(a,b)\leq \p(a+k+1,b+k)$$ for any
 $(a,b)\in\diag$  % and $(a+k+1,b+k)\in\adiag$ 
with $k=\mathrm{min}\{j\in\Z_{\geq0}\mid (a+j+1,b+j)\in\adiag\}$.
\end{lemma}
%%%%%%%%%%%%%%%%%%%%%%%%%%%%%%%%%%%%%%%%%%%%%%%%%%%%%%%

\begin{remark}
Another important generalization of plane partitions
is given by imposing 
the periodicity  $\p(u+\peri)=\p(u)$ instead of (P1).
Such plane partitions (called cylindric partitions)
are introduced and studied 
 by Gessel and Krattenthaler \cite{GK},
and they are related to the theory of quantum Schubert calculus
(see \cite{Po}).
\end{remark}
%%%%%%%%%%%%%%%%%%%%%%%%%%%%%%%%%%%%%%%%%%
\subsection{Tableaux on periodic diagrams}
%%%%%%%%%%%%%%%%%%%%%%%%%%%%%%%%%%%%%%%%%%%%5
Let $\peri\in\Z\times\F$ and $\diag\in\dset_\peri^n$.
Following \cite{SV}, we introduce 
standard tableaux on $\adiag$,
which connect the representation theory of Cherednik algebras
and the combinatorics on periodic diagrams.
\begin{definition}
\label{def;tableaux}
Let $\diag\in\dset_\peri^n$.

\smallskip\noindent
(i) A map $T:\adiag\to\Z$ is called 
a {\it tableau} on $\adiag$ if it is  a bijection
and satisfies the following condition

\smallskip
\noindent
(T1) $\ T(u+\peri)=T(u)-n\ \hbox{ for all }u\in\adiag.$

\smallskip\noindent
(ii) A tableau $T$ on $\adiag$ is called a
{\it standard tableaux} if it satisfies the following conditions:

\smallskip\noindent
(T2)
$\ T(a,b)<T(a,b+1)$ for any $(a,b)\in\adiag$ and
 $(a,b+1)\in \adiag.$

\smallskip\noindent
(T3) $\ T(a,b)<T(a+k+1,b+k)$ for any
$(a,b)\in\adiag$ and $(a+k+1,b+k)\in {\adiag}$ with
$k\in\Z_{\geq0}$.
\end{definition}
Denote by $\Tab_\peri(\adiag)$
and $\TabRC(\adiag)$
the set of tableaux %, the set of row increasing tableaux
and the set of standard tableaux on $\adiag$ respectively.
Similarly to Lemma \ref{lem;equi_PP}, we have
\begin{lemma}\label{lem;equi_TabRC}
 A tableau $T$ on $\adiag$
%map $T:\adiag\to\Z$
is a standard tableaux if and only
if it satisfies %the condition $\mathrm{(T1)}$ and
the following conditions$:$

\smallskip\noindent
$\mathrm{(T'2)}$
$\ T(a,b)<T(a,b+1)$ for any $(a,b)\in\diag$ and
 $(a,b+1)\in \diag.$

\smallskip\noindent
$\mathrm{(T'3)}$ $\ T(a,b)<T(a+k+1,b+k)$ for any
$(a,b)\in\diag$ %and $(a+k+1,b+k)\in\adiag$ 
with $k=\mathrm{min}\{j\in\Z_{\geq0}\mid (a+j+1,b+j)\in\adiag\}.$
\end{lemma}

For $T\in\Tab_\peri(\adiag)$ and $w\in\affW$, the map
$wT:{\adiag}\to\Z$ given by
$(wT)(u)=w(T(u))\ (u\in{\adiag})$
is also a tableau on ${\adiag}$.
The assignment $T\mapsto wT$ gives
an action of $\affW$ on $\Tab_\peri(\adiag)$.

\begin{proposition}[\cite{SV}]\label{pr;W-tab}
%$\mathrm{(i)}$
For any fixed $T\in\Tab_\peri(\adiag)$,
the assignment $w\mapsto wT$
gives the bijection
$\psi_T:\affW\isomto\Tab_\peri(\adiag)$ between sets.

 \end{proposition}
\begin{remark}
It is possible to give
the inverse image of $\TabRC(\adiag)$ by $\psi_T$ explicitly.
See \cite[Theorem 3.19]{SV}.
\end{remark}

We define a tableau on a (classical)
skew diagram $\fdiag$ as a  bijection
$T:\fdiag\to[1,n]$, and define  %row increasing tableaux and
standard tableaux analogously.
Denote  by $\tab(\diag)$
and $\tabRC(\diag)$ the set of tableaux
 and the set of standard tableaux on $\diag$ respectively.

A tableau $T$ on $\diag\in\dset_\peri^n$
can be extended to a tableau on $\adiag$
by the (quasi) periodicity $T(u+\peri)=T(u)-n$.

This gives an embedding
$$\tab(\diag)\to\Tab_\peri(\adiag),$$
through which we often regard $\tab(\diag)$ as a subset
of $\Tab_\peri(\adiag)$.

In the sequel, we treat the case
$\peri\in\Z_{\leq -1}\times\F$.

Then, under the identification above, $\tabRC(\diag)$
 is thought as a subset of $\TabRC(\adiag)$.
%%%%

Define a special tableau, which we call the
{\it row reading tableau} on $\diag$, % and on $\adiag$
by
\begin{equation}\label{eq;rbtab}
\rbtab (a,\mu_a+j)=\sum_{k=1}^{a-1}(\lm_k-\mu_k)+j\quad
\text{for}\ a\in[1,m],\ j\in[1,\lm_i-\mu_i],
\end{equation}
and extend it to the tableau on $\adiag$.
Here $\lm_i$ and $\mu_i$ $(i\in[1,m])$  are components of
$\lm,\mu\in\F^m$
such that $\diag=\lm/\mu$.
%as in
%the condition (D2) in Definition \ref{def;skew_diagram}.
Observe that $\rbtab\in\tabRC(\diag)
\subseteq\TabRC(\adiag)$.
\begin{example}\label{ex;rowreading}
Let $n=7$, $\peri=(-2,3)$,
and let
$\diag=\lm/\mu$ with $\lm=(5,3)$ and $\mu=(1,0)$.
Then the row reading
tableau $\rbtab$ on $\adiag$ %given above
is expressed as follows:
%
%%%%%%%%%%%%%%%%%%%%%%%%%%
\begin{center}
\setlength{\unitlength}{.06cm}
\begin{picture}(150,90)(-40,-45)
\multiput(65,40)(5,0){5}{\rule{.8pt}{.8pt}}
\put(40,20){\bx{$-\!6$}}
\put(50,20){\bx{$-\!5$}}
\put(60,20){\bx{$-\!4$}}
\put(70,20){\bx{$-\!3$}}
\put(30,10){\bx{$-\!2$}}
\put(40,10){\bx{$-\!1$}}
\put(50,10){\bx{$0$}}
\put(10,0){\bx{$1$}}
\put(20,0){\bx{$2$}}
\put(30,0){\bx{$3$}}
\put(40,0){\bx{$4$}}
\put(0,-10){\bx{$5$}}
\put(10,-10){\bx{$6$}}
\put(20,-10){\bx{$7$}}
\put(85,-3){$\;\diag$}
\put(-20,-20){\bx{$8$}}
\put(-10,-20){\bx{$9$}}
\put(0,-20){\bx{$10$}}
\put(10,-20){\bx{$11$}}
\put(-30,-30){\bx{$12$}}
\put(-20,-30){\bx{$13$}}
\put(-10,-30){\bx{$14$}}
\multiput(-35,-40)(5,0){5}{\rule{.8pt}{.8pt}}
\put(80,-3){\mbox{$\biggl. \biggr\}$}}
\linethickness{1.5pt}
\put(10,10){\line(1,0){40}}
\put(10,10){\line(0,-1){10}}
\put(0,0){\line(1,0){10}}
\put(0,0){\line(0,-1){10}}
\put(0,-10){\line(1,0){30}}
\put(30,0){\line(1,0){20}}
\put(30,0){\line(0,-1){10}}
\put(50,0){\line(0,1){10}}
%extra border
\put(30,-50){{\footnotesize Figure 2}}
\end{picture}
\end{center}
%%%%%%%%%%%%%%%%%%%%%%%%
\end{example}
%%%%%%%%%%%%%%%%%%%%%%

%%%%%%%%%%%%%%%%%%%%%%%%%%%%%%%%%%%
\subsection{Connection between standard tableaux and plane partitions}
%%%%%%%%%%%%%%%%%%%%%%%%%%%%%%%%%%%%%%%%%%%%%%%%
%\medskip
Let $\peri\in\Z_{\leq-1}\times\F$ and  $\diag\in\dset^n_\peri$.

For $T\in\Tab_\pri(\adiag)$,
define the map $\projmapop{T}:\adiag\to\Z$ by
\begin{equation}
  \label{eq;rho}
\projmapop{T}(u)=\left\lceil\frac{T(u)-1}{n}\right\rceil
\quad\hbox{for  }u\in\adiag,  
\end{equation}
 where $\lceil\, t\, \rceil$ ($t\in\Q$) denotes the maximum integer 
which is not greater than $t$.

In other words, $\projmapop{T}$ is defined as the  
map $\projmapop{T}:\adiag\to\Z$ 
such that 
$\bar T:=T-n\projmapop{T}$ is a map $\adiag\to[1,n]$
(see Figure 3).
%in $\Map(\adiag, [1,n])$.
%
%%%%%%%%%%%%%%%%%%%%%%%%%%
%\begin{figure}
\begin{center}
\setlength{\unitlength}{.06cm}
\begin{picture}(150,35)(0,-15)
\put(10,0){\bx{$0$}}
\put(20,0){\bx{$6$}}
\put(30,0){\bx{$17$}}
\put(10,-10){\bx{$8$}}
\put(20,-10){\bx{$9$}}
\put(50,-3){$=$}
\put(62,-3){$5$}
\put(70,0){\bx{$-1$}}
\put(80,0){\bx{$1$}}
\put(90,0){\bx{$3$}}
\put(70,-10){\bx{$1$}}
\put(80,-10){\bx{$1$}}
\put(108,-3){$+$}
\put(120,0){\bx{$5$}}
\put(130,0){\bx{$1$}}
\put(140,0){\bx{$2$}}
\put(120,-10){\bx{$3$}}
\put(130,-10){\bx{$4$}}
\end{picture}

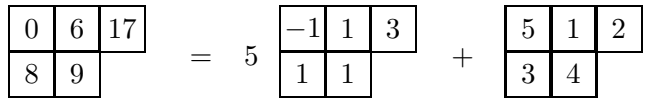
{\footnotesize Figure 3.} $T=5\projmapop{T}+\bar T$
\end{center}
%\end{figure}
%%%%%%%%%%%%%%%%%%%%%%%%
%\end{example}
%%%%%%%%%%%%%%%%%%%%%%

\smallskip
It is easy to see that if 
 $T\in\TabRC(\adiag)$ then $\projmapop{T}$ is a plane partition
on $\adiag$.
Namely, 
the assignment $T\mapsto \projmapop{T}$  gives a map 
$\projmap:\TabRC(\adiag)\to\trigPPp(\diag)$.

For $T\in\Tab_\peri(\adiag)$,
we put
\begin{align*}
 \adiag^T&=T^{-1}([1,n])\subset \adiag.
\end{align*}
The following is easy.
%%%%%%%%%%%%%%%%%%%%%%%%%%%%%%%%
\begin{lemma}
$\mathrm{(i)}$  
For any $T\in\TabRC(\adiag)$,
the set  $\adiag^T$ is a skew diagram 
consisting of $n$ elements.
 
\smallskip\noindent
$\mathrm{(ii)}$ 
For any $\p\in\trigPPp(\diag)$ and $k\in\Z$,
the set $\p^{-1}(k)$ is a skew diagram %of order $n$.
consisting of $n$ elements.

\smallskip\noindent
$\mathrm{(iii)}$ For any
 $T\in\TabRC(\adiag)$ and $k\in\Z$, it holds that
 $\projmapop{T}^{-1}(k)=\adiag^T-k\peri$.
\end{lemma}

Consider the set $\W\backslash\Tab_\peri(\adiag)$
of the $\W$-orbits on $\Tab_\peri(\adiag)$.
Let %$\Worb{\TabR(\adiag)}$ and
 $\Worb{\TabRC(\adiag)}$ %respectively 
denote
the image of %$\TabR(\adiag)$ and
 $\TabRC(\adiag)$ 
under the projection $\Tab_\pri(\adiag)\to\W\backslash\Tab_\pri(\adiag)$.
%%%%%%%%%%%%%%%%%%%%%%%%%%%%%%%%%%%%%%%%%%%%%%%%
\begin{proposition}\label{pr;TabandPP}
Let $\peri\in\Z_{\leq -1}\times\F$ and $\diag\in\dset_\peri^n$.
The map $\projmap:\TabRC(\adiag)\to\trigPPp(\diag)$
is surjective, and moreover it factors the bijection
$$\projmap:\Worb{\TabRC(\adiag)}\isomto\trigPPp(\diag).$$
\end{proposition}
%%%%%%%%%%%%%%%%%%%%%%%%%%%%%%%%%%%%%%%%%%%%%%%%%%%%%%%
\noindent{\it Proof.}
Take any $\p\in\trigPPp(\diag)$.
Then $\diag=\sqcup_{k\in\Z}\diag(k)$, where
we put $\diag(k)=
\p^{-1}(k)\cap \diag$.
Observe that
$\p^{-1}(k)$ is a skew diagram, 
%as $\p$ is row-column weakly increasing,
and so is $\diag(k)$.
One can find  a tableau $S$ on $\diag$ such that
the restriction $S|_{\diag(k)}$ is 
a standard tableau on $\diag(k)$,
%row-column increasing on $\diag(k)$
%(in the sense of (T2)(T3) in Definition \ref{def;tableaux}),
and moreover $S(u)<S(v)$ for any $u,v\in\diag$ with $\p(u)>\p(v)$.

Extend $S$ to a map $\adiag\to\Z$ periodically by $S(u+\peri)=S(u)$.
%By $\p\in\trigPPp(\diag)$, 
Then, it follows that
$T_\p:=S+n\p$ 
%is a standard tableau as $\p\in\trigPPp(\diag)$,
%and it
is a standard tableau on $\adiag$.
It is obvious that $\projmapop{T_\p}=\p$.
Therefore $\projmap$ is surjective.

Let $T\in\TabRC(\adiag)$. 
Then there exists  $\bar T\in\tab(\diag)$ such that
 $T(u)=\bar T(u)+n\projmapop{T}(u)$ for any $u\in\diag$.
%with some $\bar T\in\tab(\diag)$. %and $\p\in\trigPPp(\diag)$.
Observe that  $wT(u)=w(\bar T(u)+n\projmapop{T}(u))
=w(\bar T(u))+n\projmapop{T}(u)$  for any $w\in\W$ and $u\in\diag$. 
This implies 
that $\projmapop{T_1}=\projmapop{T_2}$ if and only if
$T_1=wT_2$ for some $w\in\W$, and hence 
the map $\projmap$ factors the bijection
$\Worb{\TabRC(\adiag)}\isomto\trigPPp(\diag)$.
\qed
%
%%%%%%%%%%%%%%%%%%%%%%%%%%%%%%%%%%%%%%%
\subsection{Content}
%%%%%%%%%%%%%%%%%%%%%%%%%%%%%%%%%%%%%%
Let $\peri\in\Z_{\leq-1}\times\F$ and
%$\peri=\mkappa$ with $m\in\Z_{\geq1}$ and $\kappa\in\F$.
%$\Q_{\leq0}$.
 $\diag\in\dset_\peri^n$.
%and $\adiag$ the corresponding periodic diagram.
Let $\con$ denote the map
$\Z\times\F \to\F$ given by $\con(a,b)=b-a$.
For a tableau $T$ on $\adiag$, 
define 
a function 
 $\con_T:\Z\to\F$
by 
$$\con_T(i)=\con(T^{-1}(i))\ \ (i\in\Z).$$

The function $\con_T$ is called the {\it content} of $T$.
%
%Put $\kappa=\con(\peri)$.
For later use, we give several lemmas below.
%, which will be used in Section \ref{sec;tabrep_trig}.
The first one is easy:
\begin{lemma}
  \label{lem;con_weq}
Let $T\in\Tab_\peri(\adiag)$.
Then

\smallskip\noindent
$\mathrm{(i)}$ $\con_T(i+n)=\con_T(i)-\con(\peri)$ for any $i\in\Z$.

\smallskip\noindent
$\mathrm{(ii)}$
 $\con_{wT}(i)=\con_T(w^{-1}(i))$
 for any $i\in\Z$ and $w\in\affW$.
\end{lemma}
\begin{lemma}
  \label{lem;con_different}
 Let $S,T\in\TabRC(\adiag)$. Then
 $\con_S= \con_T$ if and only if $S= T$.
\end{lemma}
\noindent{\it Proof.}
Let $T\in\TabRC(\adiag)$. 
Then it follows from the definition of
the standard tableaux that $T(a,b)<T(a+k,b+k)$
for any $(a,b),(a+k,b+k)\in\adiag$ with $k\in\Z_{\geq1}$.
The statement follows easily from this property.
\qed

\begin{lemma}%\cite{SV}
\label{lem;noncritical}
Let $T\in\TabRC(\adiag)$ and $i\in [0,n-1]$.

\smallskip\noindent
$\rm{(i)}$
$\con_T(i)-\con_T(i+1)\neq0$.

\smallskip\noindent
$\rm{(ii)}$
$s_iT\in\TabRC(\adiag)$ if and only if
$\con_T(i)-\con_T(i+1)\notin \{-1,1\}$.
\end{lemma}
\noindent{\it Proof.}
Follows
easily from the skew property and
the definition of the standard tableaux.
\qed
%%%%%%%%%%%%%%%%%%%%%%%%%%%%%%%%%%%%%%%%%%%
%Let $l(w)$ denote the length of $w\in\affW$.
\begin{lemma}  \label{lem;weakorder}
Let $T\in\TabRC(\adiag)$.
Let $w\in\affW$ and $i\in[0,n-1]$ such that
$l(s_iw)<l(w)$.
If $wT\in\TabRC(\adiag)$ then $s_iwT\in\TabRC(\adiag)$.
\end{lemma}
%%%%%%%%%%%%%%%%%%%%%%%%%%%%%%%%%%%%%%%%%%%%%%%%%%%%%%%%%%%%%%
\noindent{\it Proof.}
It follows from $l(s_iw)<l(w)$ that $w^{-1}(i)>w^{-1}(i+1)$.
Put $S=wT$ and suppose that
$S\in\TabRC(\adiag)$ and $s_iS\notin\TabRC(\adiag)$.
Then $\con_S(i)-\con_S(i+1)=\pm 1$.
Suppose that $\con_S(i)-\con_S(i+1)=1$.
Then it follows from $S\in\TabRC(\adiag)$ that
$S^{-1}(i+1)=S^{-1}(i)+(j+1,j)$ for some $j\in\Z_{\geq0}$.
But then we have $T(S^{-1}(i))=w^{-1}(i)>w^{-1}(i+1)=T(S^{-1}(i+1))$.
This contradicts $T\in\TabRC(\adiag)$.
Similarly we have a contradiction
in the case $\con_S(i)-\con_S(i+1)=-1$. Hence 
$s_iS=s_iwT\in\TabRC(\adiag)$.
\qed
%\smallskip

%%%%%%%%%%%%%%%%%%%%%%%%%%%%%%%
\section{Generating functions}\label{sec;generating_function}
%%%%%%%%%%%%%%%%%%%%%%%%%%%%%%%%%%%%%
Let $\diag$ be a skew diagram.
For $\p\in \PP(\diag)$, define 
\begin{align*}
%&|T|=\sum_{u\in\lm}T(u)\ \hbox{ for }T\in\Tab(\alm),\\
&|\p|=\sum_{u\in\diag}\p(u).
\end{align*}
For a subset $A$ of $\PP(\diag)$, we define the generating function 
for $S$ by
$$\GF(A;q)=\sum_{\p\in A}q^{|\p|}.$$
The purpose of this section is to compute the 
generating functions 
$\GFRperi$ and $\GFCperi$.
%%%%%%%%%%%%%%%%%%%%%%%%%%%%%%%%%%%%%%%%%%%%%%
\subsection{Count on single columns}
%%%%%%%%%%%%%%%%%%%%%%%%%%%%%%%%%%%%%%%%%%%%%%%
%
%\def\r{p}
\def\Pset{M}
As a first example,
we compute the generating functions
for   %the tableaux representation corresponding to 
a special skew diagram
$$\vert=\{(a,1)\in\Z\times\Z \mid a\in[1,n]\}$$ 
by a naive enumerative method. 
We will see in Section \ref{ss;rep_single_column}
that they are related with character formulas for
finite-dimensional representations of the rational Cherednik
algebra of type $A$ obtained by Berest-Etingof-Ginzburg \cite{BEG}.

It turns out that
the computation is easy for most 
$\kappa\in\Qgeq$  (see Lemma \ref{lem;unless_n/r}), 
and
the only interesting case is
\begin{equation}
  \label{eq;kappa=n/r}
\kappa=n/\r\hbox{ with }
\r\in\Z_{\geq1},\ (n,\r)=1,
\end{equation}
which we treat in the sequel.
%%%%%%%%%%%%%%%%%%%%%%%%%%%%%%%%%%%%%%%%%%%%%%%%%%%%%%%%%%%%%%%%
For $k\in\Z_{\geq1}$, put
\begin{equation}
  \label{eq;qfactorial}
  \qfac{k}=(1-q)(1-q^2)\dots (1-q^k),
\end{equation}
and put $\qfac{0}=1$.
%%%%%%%%%%%%%%%%%%%%%%%%%%%%%%%%%%%%%%%%%%%%%%%%%%%%%%
\begin{proposition}\label{pr;GF_single_column}
Let $\r\in\Z_{\geq1}$ such that $(n,r)=1$.
Put $\peri=(-n,n/\r-n)$.
Then 
\begin{align*}
\GF(\PPp(\vert);q)=\GF(\PPRp(\vert);q)
%\sum_{\p\in\PPRp(\vert)_{\geq0}}q^{|\p|}
&=
\frac{ \qfac{n+\r-1} }{ \qfac{n} \qfac{\r} },
\end{align*}
\begin{align*}
\GF(\PPCp(\vert);q)
%\sum_{\p\in\PPCp(\diag)_{\geq0}}q^{|\p|}
&=
\begin{cases}
\frac{ \qfac{\r-1}}{ \qfac{n} \qfac{\r-n} }\ 
\ &\hbox{ if }\r\geq n,\\
0\ &\hbox{ if }\r<n.
\end{cases}
\end{align*}
\end{proposition}
%%%%%%%%%%%%%%%%%%%%%%%%%%%%%%%%%%%%%%%%%%%%%%%%%%%%%%%%%%%%%%%%%
\noindent{\it Proof.}
Let $\p\in\PP(\vert)$.
It follows from Lemma \ref{lem;equi_PP} that
$\p\in\PPp(\vert)$ if and only if 
$\p(n,1)\leq \p(1+rn,1)=\p(1,1)+r$.
Identifying the diagram $\vert$ with the set $[1,n]$ 
via the correspondence $(a,1)\mapsto a$, we have
\begin{align*}
\PPp(\vert) %_{\geq0}
=\PPRp(\vert)
&=\{\p:[1,n]\to \Z_{\geq0}\mid
\p(1)\leq \p(2)\leq \dots \leq \p(n)\leq \p(1)+\r\},\\
\PPCp(\vert) %_{\geq0}
&=\{\p:[1,n]\to \Z_{\geq0}\mid
\p(1)<\p(2)<\dots <\p(n)<\p(1)+\r\}.
\end{align*}
Note that
$\PPRp(\vert)=\sqcup_{k\in\Z_{\geq0}}
\PPRp(\vert)_{k}$ and
$\PPCp(\vert)=\sqcup_{k\in\Z_{\geq0}}
\PPCp(\vert)_{k}$,
where
\begin{align*}
%\end{align*}
%\begin{align*}
\PPRp(\vert)_{k}
&=\{\p\in\PPRp(\vert)\mid \p(1)=k\}  \\
&\cong \{\p:[2,n]\to [k,k+\r]\mid \p(2)\leq\dots \leq\p(n)\},\\
\PPCp(\vert)_{k}
&=\{\p\in\PPCp(\vert)\mid \p(1)=k\}  \\
&\cong \{\p:[2,n]\to [k+1,k+\r-1]\mid \p(2)<\dots <\p(n)\}.
\end{align*}
The generating functions for these sets
have been calculated classically (see e.g., \cite[Section 1.3]{St1}).
In particular for $k=0$, we have
\begin{align*}
%\frac{\qfac{r-1}}{\qfac{n-1}\qfac{r-n}}.
&\GF(\PPRp(\vert)_{0};q)
=\frac{\qfac{n+\r-1}}{\qfac{n-1}\qfac{\r}},\\
&\GF(\PPCp(\vert)_{0};q)
%\sum_{\p\in\PPCp(\vert)_{0}}q^{|\p|}=
=
\begin{cases}
\frac{ \qfac{\r-1}}{ \qfac{n-1} \qfac{\r-n}}\ 
\ &\hbox{ if }\r\geq n,\\
0\ &\hbox{ if }\r<n.
\end{cases}
\end{align*}
%
%Obviously, we have
Noting
the equalities
 $\GF(\PPRp(\vert)_{k};q)=q^{kn}\GF(\PPRp(\vert)_{0};q)$
and 
$\GF(\PPCp(\vert)_{k};q)=q^{kn}\GF(\PPCp(\vert)_{0};q)$,
the statement follows.
\qed
%%
%%%%%%%%%%%%%%%%%%%%%%%%%%%%%%%%%%%%%%%%%%%%%%%%%%%%%%%%%%%%%%%%
\subsection{Kostka polynomial}
%%%%%%%%%%%%%%%%%%%%%%%%%%%%%%%%%%%%%%%%%%%%%%%%%%
\def\ccwt{{\mathbf{h}}}
\def\cocharge{{\mathbf{cc}}}

The main purpose in the rest of Section \ref{sec;generating_function}
is to compute the generating function
\begin{align*}
&\GFCperi=\sum_{\p\in\PPCp(\diag)}
q^{|\p|}
%&
\end{align*}
%for $\PPCp(\diag)$
when $\adiag\subset\Z\times\Z$.
It turns out that $\GFCperi$ is expressed by level restricted
Kostka polynomials,
and it is proved using 
Lascoux-Sch\"utzenberger 
type expression \cite{LS},
which we will see  below.

\def\cccoef{d}
Put
\begin{align*}
&\dsetint^n=\dset^n\cap (\Z\times\Z),\ \ \dsetint_\peri^n=\dset_\peri^n\cap (\Z\times\Z).
\end{align*}
Let $\diag\in\dsetint^n$.
For $T\in\tab(\diag)$
 and $i\in[1,n]$,
define
 \begin{equation}\label{eq;cocharge}
 \cccoef_i(T)=\begin{cases}
 &1\quad %\hbox{if}\ a<a'
\ \hbox{ if }
T^{-1}(i+1)-T^{-1}(i)\in\Z_{\geq1}\times \Z,
\\
 &0\quad %\hbox{if}\ a\geq a'
 \ \hbox{ otherwise,}
 \end{cases}
%\ \hbox{ where }(a,b)=T^{-1}(i),\ (a',b')=T^{-1}(i+1),
 \end{equation}
 and put
\begin{equation}
  \label{eq;eta_T}
\ccwt_T=\sum_{i\in[1,n]}
\left(\sum_{k\in[1,i-1]}\cccoef_k(T)\right)\e_i\in P\;(=\oplus_{i\in[1,n]}
\Z\e_i).
\end{equation}

Then  $\brac{\ccwt_T}{\alch_i}=-\cccoef_{i}(T)$ 
for $i\in[1,n]$ and hence $\ccwt_T\in P^-$, where
we put
\begin{align*}
P^-&=\{\zeta\in P\mid \brac{\zeta}{\alch}\geq0
\hbox{ for all } \al\in R^+\}.
\end{align*}

For $\zeta\in P$, 
we write $|\zeta|=\sum_{i\in[1,n]}\brac{\zeta}{\ech_i}$.
%%%%%%%%%%%%%%%%%%%%%%%%%%%%%%%

%\medskip
Define
\begin{align} \label{eq;LS1}
%&{K}_{\lm\; (1^n)}(q)=
%\sum_{T\in\tabRC(\dlm)}q^{\charge(T)},\\
&{\check K}_{\diag\; (1^n)}(q)=
\sum_{T\in\tabRC(\diag)}q^{|\ccwt_T|},\\
 \label{eq;LS2}
&{K}_{\diag\; (1^n)}(q)=q^{\frac{1}{2} n(n-1)}
 {\check K}_{\diag\; (1^n)}(q^{-1}),
%\check{K}_{\diag'\, (1^n)}(q),
\end{align}
where $(1^n)$ denote the partition $(1,1,\dots,1)$ of $n$,
The polynomial $K_{\diag\, (1^n)}(q)$
is called the Kostka polynomial
associated with the skew diagram $\diag$ and 
the partition $(1^n)$.
 It is known that
\begin{equation} \label{eq;LSconj}
K_{\diag\; (1^n)}(q)=\check{K}_{\diag'\, (1^n)}(q), 
\end{equation}
where $\diag'$ is the conjugate of $\diag$:
$\ \diag'=\{(a,b)\mid (b,a)\in\diag\}$.

%%%%%%%%%%%%%%%%%%%%%%%%%%%%%%%%%%%%%%%%%%%%%%%%%%%%
\begin{remark}
 For partitions $\lm$ and $\mu$, 
the Kostka polynomials $K_{\lm\,\mu}(q)$ can be defined 
as the coefficients in the expansion
$$s_\lm(z_1,\dots,z_n)=\sum_{\mu\vdash n}K_{\lm\,\mu}(q)
p_\mu(z_1,\dots,z_n;q),$$
where $s_\lm$ is the Schur polynomial and 
$p_\mu$ is the Hall-Littlewood polynomial
$($see e.g., \cite{Mac}$)$.
The expression of the Kostka polynomials of the form 
\eqref{eq;LS1}\eqref{eq;LS2} is 
given in \cite{LS}  
\end{remark}
%%%%%%%%%%%%%%%%%%%%%%%%%%%%%%%%%%%%%%%%%%%%%%%%%%%%%%%%%%%%
\subsection{Generating functions for classical plane partitions}
%%%%%%%%%%%%%%%%%%%%%%%%%%%%%%%%%%%%%%%%%%%%%%%%%%%%%%%%%%%%%%%%
Concerning the computation of the generating function $\GFC$
% the set  $\PPC(\diag)$ of
for (classical) column strict plane partitions,
several approaches have been known.
%the generating function $\GFC$.
(See \cite{St;pp} \cite[Section 7.21]{St2} \cite[Example 1.5.12]{Mac}.)
We show the formula for  $\GFC$
using the following statement,
and will generalize it to the periodic case later.

Recall that we identify
$\tabRC(\diag)$ 
as a subset of $\TabRC(\adiag)\subset \Map(\adiag,\Z)$.
%%%%%%%%%%%%%%%%%%%%%%%%%%%%%%%%%%%%%%%%%%%%%%%%%%%%%%%%%%%%
\begin{proposition}\label{pr;PtabRC_PPC}
Let $\diag\in\dsetint^n$.
  The assignment $(\zeta,T)\mapsto \projmapop{t_{\zeta+\ccwt_T} T}$,
where $\projmap$ is given by \eqref{eq;rho},
gives a bijection
\begin{align*}
  &P^-\times\tabRC(\diag)\isomto \trigPPC(\diag).
\end{align*}
\end{proposition}
%%%%%%%%%%%%%%%%%%%%%%%%%%%%%%%%%%%%%%%%%%%%%%%%%%%%%%%%%%%%%%
\noindent{\it Proof.}
It follows from the definition of $\ccwt_T$ 
that $\projmapop{t_{\ccwt_T}T}\in \trigPPC(\diag)$
for any $T\in\tabRC(\diag)$.
Hence 
$\projmapop{t_{\zeta+\ccwt_T}T}\in \trigPPC(\diag)$ for any 
$\zeta\in P^-$ and $T\in\tabRC(\diag)$.

To show that this map is a bijection,
we define a map $\beta:\trigPPC(\diag)\to P^-\times \tabRC(\diag)$
as follows.
Take $\p\in\trigPPC(\diag)$.
Introduce an order in the set $\diag$ by
\begin{equation}\label{eq;order_diagram}
(a,b)\prec (a',b')
\Leftrightarrow 
\begin{cases}
\p(a,b)<\p(a',b'),\hbox{ or }\\
\p(a,b)=\p(a',b')\hbox{ and }b<b'. 
\end{cases}
\end{equation}
Observe that this gives a total order because $\p$ is column strict.

Let $u_1,\dots,u_n$ be the elements in $\diag $
such that $u_1\prec u_2\prec\dots \prec u_n$.
%$i<j\Leftrightarrow \p(u_i)\prec \p(u_j)$. 
Define $T_\p\in\tab(\diag)$ by $T_\p(u_i)=i$,
and  define $\zeta_\p\in P$
by $\brac{\zeta_\p}{\ech_i}=\p(u_i)-\brac{\ccwt_{T_\p}}{\ech_i}$.
Then it is easy to check that $T_\p\in\tabRC(\diag)$.

Let us check $\zeta_\p\in P^-$.
We have 
$$\brac{\zeta_\p}{\alch_{i}}=\p(u_i)-\p(u_{i+1})-
\brac{\ccwt_{T_\p}}{\alch_i}.$$
Observe that $\brac{\ccwt_{T_\p}}{\alch_i}=-\cccoef_i(T_\p)\geq-1$.

Suppose that $\brac{\zeta_\p}{\alch_{i}}>0$.
This occurs only if $\p(u_i)-\p(u_{i+1})=0$ and 
$\brac{\ccwt_{T_\p}}{\alch_i}=-1$.
The latter equality implies
 $u_{i+1}-u_i\in\Z_{\geq1}\times\Z$ by \eqref{eq;cocharge}.
%%and  
%$\brac{\ccwt_{T_\p}}{\alch_i}=-1$ 
By $\p(u_i)=\p(u_{i+1})$, 
it follows from the definition \eqref{eq;order_diagram}
of the order $\prec$ that
$u_{i+1}-u_{i}\in\Z\times \Z_{\geq1}$.
Hence $u_{i+1}-u_{i}\in\Z_{\geq1}\times \Z_{\geq1}$,
and this implies $\p(u_{i+1})>\p((u_i)$ 
as $\p$ is column strict.
This is a contradiction.
Therefore %$\p(u_i)<\p(u_{i+1})$ and hence
$\brac{\zeta_\p}{\alch_{i}}\geq0$ for any $i\in[1,n]$,
namely, $\zeta_\p\in P^-$.

We set $\beta(\p)=(\zeta_\p,T_\p)$.
%It is easy to see that
It is easy to see that $\beta$ gives an inverse map of the map
$(\zeta,T)\mapsto \projmapop{t_{\zeta+\ccwt_T}T}$.
\qed

\medskip
By restricting the map in Proposition \ref{pr;PtabRC_PPC}
we have a bijection
\begin{equation}\label{eq;PregtabRC_PPC}
\Preg^-\times\tabRC(\diag)\isomto \PPC(\diag), %_{\geq0},
\end{equation}
where we put
$$\Preg^-
=\{\zeta\in P^-\mid \brac{\zeta}{\ech_1}\geq0\}.$$
As a consequence,
%%%%%%%%%%%%%%%%%%%%%%%%%%%%%%%%%%%%%%%%%%%%%%%%%%%%%%%%%%%%%%%%%%%%%%%%
we have the following well-known 
result. %(See e.g., \cite[Examples 1.5.12]{Mac} \cite[Section 7.21]{St2}).
\begin{proposition}[\cite{St;pp}]
\label{pr;char_classical}
  %Let $\kappa\in\F\setminus\{0\}$, $m\in\Z_{\geq1}$ and
Let $\diag\in\dsetint^{n}$.
Then
\begin{align*}
\GFR&=
\frac{1}{\qfac{n}}
{K}_{\diag\; (1^n)}(q),\\
\GFC&=
\frac{1}{\qfac{n}}
{\check K}_{\diag\; (1^n)}(q).
\end{align*}
\end{proposition}
%%%%%%%%%%%%%%%%%%%%%%%%%%%%%%%%%%%%%%%%%%%%%%%%%%%%%%%%%%%%%%%%
%
\noindent{\it Proof.}
Note that $|\projmapop{t_{\zeta+\ccwt_T}T}|=
\sum_{i\in[1,n]}\brac{\zeta+\ccwt_T}{\ech_i}
%=|\zeta+\ccwt_T|
=|\zeta|+|\ccwt_T|$.
The bijection \eqref{eq;PregtabRC_PPC} derived from
 Proposition \ref{pr;PtabRC_PPC}
 implies 
%the bijection
% $\Preg^-\times\tabRC(\diag)\isomto \PPC(\diag)_{\geq0}.$
$$\GFC=\sum_{\zeta\in \Preg^-}\sum_{T\in\tabRC(\diag)}
q^{|\zeta|+|\ccwt_T|}=
\frac{1}{\qfac{n}}
{\check K}_{\lm\; (1^n)}(q).$$

%where $\lm'$ is the conjugate of $\lm$.
To show the formula for $\GFR$,
we consider the conjugate operation.
For $\p\in\PPR(\diag)$, 
the function $\p'(a,b)=\p(b,a)$ for $(b,a)\in\diag$
gives a column strict plane partition on $\diag'$.
We have $|\p|=|\p'|$ and 
the correspondence $\p\mapsto \p'$ gives a bijection
$\PPR(\diag)\isomto \PPC(\diag')$.
Hence $\GF(\PPR(\diag);q)=\GF(\PPC({\diag'});q)
=\check K_{\diag'\, (1^n)}(q)=K_{\diag\,(1^n)}(q)$.
\qed

%%%%%%%%%%%%%%%%%%%%%%%%%%%%%%%%%%%%%%%%%
\subsection{Restricted tableaux}
\label{ss;restricted_tableaux}
%%%%%%%%%%%%%%%%%%%%%%%%%%%%%%%%%%%%%%%%%%%%%%%%%%
%%%
\def\dom{{{}_\Z{\mathfrak D}}}
In the rest of this section 
we restrict ourselves to the case where $\diag\in\dsetint_\peri^n$
with $\con(\peri)\in\Z_{\geq1}.$

Let $\kappa\in\Z_{\geq1}$ and $m\in\Z_{\geq1}$.
Put $\peri=\mkappa$ and let $\diag\in\dsetint^{n}_\peri$.
Then $\diag$ is expressed as $\diag=\lm/\mu$
for some  $\lm,\mu\in\domset_{m,\kappa}$.
Here, recall that $\domset_{m,\kappa}$ is defined 
as the subset of $\F^m$ consisting of 
the elements satisfying the dominance condition
\eqref{eq;domsetmk1}\eqref{eq;domsetmk2}.

For $T\in\tabRC(\diag)$, we have a sequence 
$T^{-1}(\{1\}),\ T^{-1}(\{1,2\}),\dots,
T^{-1}([1,n])$
%$T^{-1}([1,k])$ $(k=1,2,\dots,n)$ 
of skew diagrams.
Define $\lm_T^{(k)}\in\F^m$ as  the unique element such that
$T^{-1}([1,k])=\lm_T^{(k)}/\mu$.

\begin{definition}
A standard tableaux $T$ on $\diag$ is said to be
{\it $\peri$-restricted} %(or simply, {\it restricted})
if $\lm_T^{(k)}\in \domset_{m,\kappa}$
for all $k\in [1,n]$.
\end{definition}
Define
\begin{align}
  \restab_\pri(\diag)&=\{T\in\st(\diag)\mid T
\hbox{ is } \peri\hbox{-restricted}\}.
%,\\
\end{align}
%
%%%%%%%%%%%%%%%%%%%%%%%%%%%%%%%%%%%%%%
\begin{lemma}\label{lem;leq2}
$\mathrm{(i)}$ If $T\in\restab_\pri(\diag)$,
then for any $i,j\in[1,n]$ with $i>j$, it holds that
 $\con_T(j)-\con_T(i)\leq \kappa-2$.

\smallskip\noindent
$\mathrm{(ii)}$  If $T\in\tabRC(\diag)$ and $T\notin\restab_\pri(\diag)$,
then there exist $i,j\in[1,n]$ such that $i>j$
and $\con_T(j)-\con_T(i)=\kappa-1$.
\end{lemma}
%%%%%%%%%%%%%%%%%%%%%%%%%%%%%%%%%%%%%%%%%%%%%
\noindent
{\it Proof.}
(i) Suppose that $T\in\restab_\pri(\diag)$.
Take any $i\in[1,n]$.
It follows from $\lm_T^{(i)}\in \domset_{m,\kappa}$ that 
$\con_T(j)-\con_T(i)\leq\kappa-1$ for any $j\in[1,i-1]$.
Suppose that the equality holds.
Then, it is easy to see that $T^{-1}(j)-T^{-1}(i)\in\Z_{\leq-1}\times \Z$
(namely, $T^{-1}(i)$ is located below $T^{-1}(j)$).
%Then 
%$T^{-1}(i)=(m,b)$ and $T^{-1}(j)=(1,b+\kappa-m)$ for some $b\in\Z$.
But this implies $\lm_T^{(i-1)}\notin \domset_{m,\kappa}$,
and this is a contradiction.

\noindent
(ii) Suppose that $T\notin\restab_\pri(\diag)$.
%Since $T^{-1}([1,n])$ is restricted,
Since $\lm_T^{(n)}=\lm\in\domset_{m,\kappa}$, 
there exists $i\in[1,n]$ such that
$\lm_T^{(i-1)}\notin\domset_{m,\kappa}$
and $\lm_T^{(i)}\in\domset_{m,\kappa}$.
It is easy to see that %there exists $j\in[1,i-1]$ such that
$\con_T(j)-\con_T(i)=\kappa-1$,
where  $j$ is a number such that 
$\con_T(j)=\mathrm{max}\{\con_T(k)\mid k\in[1,i-1]\}$. 
\qed

\medskip
Recall that we regard $\tabRC(\diag)$ as a subset 
of $\TabRC(\adiag)$.
%%%%%%%%%%%%%%%%%%%%%%%%%%%%%%%%%%%%%%%%%%%%%
\begin{proposition}\label{pr;tabres_PPC}
Let $m\in\Z_{\geq1}$, $\kappa\in\Z_{\geq1}$ and
 $\diag\in\dsetint_\peri^{n}$ with $\peri={\mkappa}$.

\smallskip\noindent
$\mathrm{(i)}$
The assignment $(\zeta,T)\mapsto t_{\zeta+\ccwt_T} T$
gives an embedding
$$P^-\times \restab_\pri(\diag)\to\TabRC(\adiag).$$
$\mathrm{(ii)}$
The assignment  $(\zeta,T)
%\mapsto t_{\zeta+\ccwt_T} T
\mapsto
\projmapop{t_{\zeta+\ccwt_T} T}$
gives  bijections
$$
P^-\times \restab_\pri(\diag)\isomto \trigPPCp(\diag),
\quad
\Preg^-\times \restab_\pri(\diag)\isomto \PPCp(\diag).
$$
\end{proposition}
%%%%%%%%%%%%%%%%%%%%%%%%%%%%%%%%%%%%%%%%%%%%%%%%%%%%%%%%%
%%%%%%%%%%%%%%%%%%%%%%%%%%%%%%%%%%%%%%%%%%%%%
\noindent{\it Proof.}
In this proof, we put 
$\xi_i=\bra\zeta+\ccwt_T\mid \ech_i\ket$ for $i\in[1,n]$.
%and $s'=\bra\zeta+\ccwt_T\mid \ech_j\ket$.

\smallskip\noindent
%{\underline{Proof of (a)$\Rightarrow$(b).}}
(i)
First we shall prove 
that the image of $P^-\times\restab_\pri(\diag)$
is included in $\TabRC(\adiag)$.

Let $\zeta\in P^-$ and $T\in\tabRC(\diag)$.
Put $S=t_{\zeta+\ccwt_T} T$.

Let us check the row increasing condition (T2) 
in Definition \ref{def;tableaux}.
Take any $(a,b),(a,b+1)\in\diag$.
Put $i=T(a,b)$ and $j=T(a,b+1)$.
Then we have  $i<j$ as $T\in\tabRC(\diag)$,
and we have $\xi_i<\xi_j$ as  $\zeta+\ccwt_T\in P^-$.
Hence 
\begin{align*}
 S(a,b)&=t_{\zeta+\ccwt_T}T(a,b)=i+n\xi_i
<j+n\xi_j=S(a,b+1).
\end{align*}
Therefore  $S$ satisfies the condition (T2).

Let us check the condition (T3).
Since $\adiag\subset\Z\times\Z$,
it is enough to check
% the following condition:
%(T3'') in Remark \ref{rem;tab_integral}, namely,
%\smallskip
 $S(u)<S(v)$ for any  $u=(a,b)\in\diag$ and
 $v=(a+k+1,b+k)\in\adiag$
with $k\in\Z_{\geq0}$.
%\smallskip

Put $i=T(u)\in[1,n]$.
%Note that $v$ or $v+\peri$ is in $\diag$.
Let $p\in\Z_{\geq0}$ be the number such that
$v+p\peri\in\diag$.

Suppose that $p=0$, namely,
$v\in\diag$.
Then $T(v)=j$ for some
$j\in[1,n]$, and, by similar argument as above,
we have  $i<j$ and  $\xi_i\leq \xi_j$,
and hence $S(u)<S(v)$.

Suppose that $v+p\peri\in\diag$ with $p>0$.
Then $T(v)=j+pn$ for some $j\in[1,n]$.
 %using $\con_T(j)=\con_T(j+n)+\kappa$,
We have
$\con_T(j)-\con_T(i)=\con_T(j+pn)+p\kappa-\con_T(i)=p\kappa-1>\kappa-2$.
By $T\in\restab_\pri(\diag)$, 
Lemma \ref{lem;leq2} implies $i<j$, and hence we have
 $\xi_i\leq\xi_j$ as $\zeta+\ccwt_T\in P^-$.
We have
\begin{align*}
 S(u)&=t_{\zeta+\ccwt_T}T(u)=i+n\xi_i<j+n(\xi_j+p)
=t_{\zeta+\ccwt_T}T(v)=S(v).
\end{align*}
Hence $S$ satisfies the condition (T3).
We have proved $S\in\TabRC(\adiag)$.

Now, we shall prove the statement (ii).
Note that this implies
the injectivity of the map 
$P^-\times \restab_\pri(\diag)\to\TabRC(\adiag)$.

\smallskip\noindent
(ii)
%\noindent
%\underline{Proof of .}
% Suppose (b). 
({\it Step 1})
We take $\zeta\in P^-$ and $T\in\restab_\pri(\diag)$,
and shall prove that $\projmapop{t_{\zeta+\ccwt_T} T}\in
\trigPPCp(\diag)$.

Put $S=t_{\zeta+\ccwt_T} T$.
We have shown $S\in\TabRC(\adiag)$, and hence
$\projmap(S)\in\trigPPp(\diag)$ by Proposition \ref{pr;TabandPP}.
Assume that $\projmapop{S}\notin\trigPPCp(\diag)$.
Then there exist 
$u=(a,b)\in\adiag$ and $v=(a+k+1,b+k)\in\adiag$ with $k\in\Z_{\geq0}$
such that $\projmapop{S}(u)=\projmapop{S}(v)$.

We may assume without loss of generality that $\projmapop{S}(u)=0$.
Then, putting $i=S(u)$ and $j=S(v)$, we have $i,j\in[1,n]$.
We have
 $T(u)=t_{-\zeta-\ccwt_T}S(u)=i-n\xi_i$
and, similarly, $T(v)=j-n\xi_j$.

Since $T\in\st(\diag)\subset \TabRC(\adiag)$,
it must hold that
$\xi_i\geq \xi_j$.
On the other hand, since we have shown that $S\in\TabRC(\adiag)$,
it must hold that $i<j$,
and hence $\xi_i\leq \xi_j$ as $\zeta+\ccwt_T\in P^-$.
Therefore we have $\xi_i=\xi_j$.

Look at two elements
 $T^{-1}(i)=u-\xi_i\peri$ 
and $T^{-1}(j)=v-\xi_j\peri$ in $\diag$.
Since $i<j$ and $a+\xi_i m< a+k+1+\xi_i m
=a+k+1+\xi_j m$, 
it follows from the definition of $\ccwt_T$ that
 $\brac{\ccwt_T}{\ech_i}<
\brac{\ccwt_T}{\ech_j}$.
Combining with $\zeta\in P^-$, this implies $\xi_i<\xi_j$.
This is a contradiction.
Therefore $\projmapop{S}\in\trigPPCp(\diag)$.

\smallskip\noindent
({\it Step 2})
Recall the bijection
$P^-\times \tabRC(\diag)\isomto \trigPPC(\diag)$
in Proposition \ref{pr;PtabRC_PPC}.
What we have to show is the surjectivity of the map
$P^-\times \restab_\pri(\diag)\to \trigPPCp(\diag)$.

Let $T\in\tabRC(\diag)$ and $\zeta\in P^-$,
and suppose that $\projmapop{t_{\zeta+\ccwt_T} T}\in\trigPPCp(\diag)$.
We shall show that $T\in\restab_\pri(\diag)$.
Note that this will complete the proof of the statement. 

Suppose that 
$T\notin \restab_{\peri}(\diag)$.
Then, 
by Lemma \ref{lem;leq2}, there exist $i,j\in[1,n]$  such that 
%$T^{-1}(i)=(m,b)$ for some $b$ and
$i>j$ and $\con_T(j)-\con_T(i)=\kappa-1$.
Combining with $T\in\tabRC(\diag)$, we have  $a>a'$,
where
we put $T^{-1}(i)=(a,b)$ and $T^{-1}(j)=(a',b')$.

Observe that 
$(a',b')-\peri=(a,b)+(k+1,k)$,
where
$k=m-1-(a-a')\in\Z_{\geq0}$.

By $i>j$ and $a>a'$, 
it follows from the definition of $\ccwt_T$ that
 $\bra\ccwt_T\mid \ech_i\ket>\bra\ccwt_T\mid \ech_j\ket$,
and hence $\xi_i>\xi_j$.

We have $t_{\zeta+\ccwt_T} T(a,b)=i+n\xi_i$ and
\begin{align*}
t_{\zeta+\ccwt_T} T(a+k+1,b+k)&=t_{\zeta+\ccwt_T} T((a',b')-\peri)\\
&=t_{\zeta+\ccwt_T}(j+n)=j+n(\xi_j+1),
\end{align*}
and hence
%\begin{align*}
%&
$\projmapop{t_{\zeta+\ccwt_T} T}(a,b)=\xi_i$ and 
$\projmapop{t_{\zeta+\ccwt_T} T}(a+k+1,b+k)=\xi_j+1$.
%\end{align*}
The assumption $\projmapop{t_{\zeta+\ccwt_T} T}\in\trigPPCp(\diag)$
implies $\xi_i<\xi_j+1$.
This is a contradiction and hence
$T\in\restab_\pri(\diag)$.
The statement has been proved.
\qed

%%%%%%%%%%%%%%%%%%%%%%%%%%%%%%%%%%%%%%%%%%%%%%%%%%%%%%%%%%%%%%%%%%%%%%%%%%%%
\subsection{Generating functions and level restricted Kostka polynomials}
%%%%%%%%%%%%%%%%%%%%%%%%%%%%%%%%%%%%%%%%%%%%%%%%%%%%%%%%%%%%
%\medskip
%Let $\kappa\in\Z_{\geq1}$.
Define
\begin{align} \label{eq;LSres}
&{\check K}_{\diag\; (1^n)}^{(\kappa-m)}(q)=
\sum_{T\in\restab_{\peri}(\diag)}q^{|\ccwt_T|}.
\end{align}
The polynomials given by
$
K_{\diag\; (1^n)}^{(\kappa-m)}(q)=
q^{\frac{1}{2} n(n-1)}
 {\check K}_{\diag\; (1^n)}^{(\kappa-m)}(q^{-1})
$
is called the
level restricted Kostka polynomial of level $\kappa-m$
associated with the skew diagram $\diag$ and the partition $(1^n)$ 
(see e.g. \cite{SS}).
We obtain a periodic analogue (or level restricted analogue) of 
the classical formulas in Proposition \ref{pr;char_classical}.
%%%%%%%%%%%%%%%%%%%%%%%%%%%%%%%%%%%%%%%%%%%%%%%%%%%%%%%%%%%%%%%
\begin{theorem}\label{th;GFCperi}
  Let $\kappa,\ m\in\Z_{\geq1}$ and put $\peri=\mkappa$.
Let $\diag\in\dsetint_\peri^{n}$.
Then
$$\GFCperi=
\frac{1}{\qfac{n}}
{\check K}^{(\kappa-m)}_{\diag\; (1^n)}(q),$$
\end{theorem}
%%%%%%%%%%%%%%%%%%%%%%%%%%%%%%%%%%%%%%%%%%%%%%%%%%%%%%%%%%%%%%%%%%%
%
\noindent{\it Proof.}
Follows from Proposition \ref{pr;tabres_PPC}-(ii) using
$|\projmapop{t_{\zeta+\ccwt_T}T}|=|\zeta+\ccwt_T|=|\zeta|+|\ccwt_T|$.
\qed
\begin{remark}\label{rem;row}
The formula for $\GFRperi$ (for $\kappa\in\Z_{\geq1}$)
has not been obtained by the argument above.
Classically, to compute the generating function $\GFR$
for the set of row strict plane partitions 
is an equivalent problem
with to compute
$\GFC$ for the set of column strict plane partitions,
since they are transferred to each other through the conjugate (transpose) 
operation.

On the other hand, for plane partitions on a periodic diagram, 
we do not have $\PPRp(\diag)\cong \PPCp(\diag)$ 
in general any more.
Actually, plane partitions (and standard tableaux) on 
a periodic diagram of period $\peri$
are transformed by the conjugation into 
those on the conjugated periodic diagram of period 
$\peri'$, for which 
we have $\con(\peri')=-\con(\peri)=-\kappa$.
For a period with negative content, our method 
can be applied to compute $\GF(\PPR_{\peri'}(\diag');q)$,
but fails to compute $\GF(\PPC_{\peri'}(\diag');q)$.\footnote{
Another approach to compute the generating functions
is to use an interpretation of plane partitions in terms of
 non-intersecting lattice paths. 
%\cite{GK}.
By this method, Gessel and Krattenthaler
obtained a determinant expression 
of the generating functions for cylindric partitions
(\cite[Theorem 2, Theorem 3]{GK}).
Connection between their results
and our formula in Theorem \ref{th;GFCperi}
will be discussed in another place.}
\end{remark}
%%%%%%%%%%%%%%%%%%%%%%%%%%%%%%%%%%%%%%%%%%%%%%%%%%%%%%%%
\section{Tableaux representations of the trigonometric Cherednik algebra}
\label{sec;tabrep_trig}
%%%%%%%%%%%%%%%%%%%%%%%%%%%%%%%

We apply combinatorics on periodic diagrams
%studied in the previous part of this article
to the study of the representation theory of Cherednik algebras.

First, we give a combinatorial construction
of a class of representations 
of the trigonometric Cherednik algebra
by modifying the construction in \cite{SV}
for the double affine Hecke algebra.

%%%%%%%%%%%%%%%%%%%%%%%%%%%%%%%%%%%%%%%%%%%%%%%%%%%%%%%%%%%%%%%%%%%%%%%%
\subsection{Trigonometric Cherednik algebra}
%%%%%%%%%%%%%%%%%%%%%%%%%%%%%%%%%%%%%%
%%%%%%%%%%%%%%%%%%%%%%%%%%%%%%%%%%%%%%%%%%%%%%%
%
We introduce several more notations.

For a set $G$, we let 
$\F G$ denote the vector space 
of $\F$-valued functions on $G$. %, where $G$ is a set.
In particular, if $G$ is a (semi)group,  then $\F G$ 
denotes the group algebra.
%%%%%%%%%%%%%%%%%%%%%%%

We let $\X$ denote the polynomial ring $\F[x_1,\dots,x_n]$,
and let $\XX$ denote 
the Laurent polynomial ring $\F[x_1^{\pm1},\dots,x_n^{\pm1}]$.

For a vector space $V$, we let $S(V)$ denote the symmetric algebra of $V$.

For an algebra $A$, we denote by $A$-$\mod$
the category of finitely generated $A$-modules.

For an element $\eta=\sum_{i\in[1,n]}\eta_i \e_i$ of $P$, we denote by 
$x^\eta=x_1^{\eta_1}x_2^{\eta_2}\dots x_n^{\eta_n}$
the corresponding element in the group algebra $\F P$.
Via the correspondence
$\eta\mapsto x^\eta$, we often identify $\F P$ with
the Laurent polynomial ring $\XX$.
The action of $\W$ on  
$\F P=\XX$ induced from the action on $P$ is
given by the permutation of variables $x_1,\dots,x_n$.

\begin{definition}
For $\kappa\in \F$,
the {\it trigonometric Cherednik algebra} 
(or the {\it degenerate double affine Hecke algebra})
$\widetilde\H_\kappa$ of type $GL_n$
is defined as the unital associative $\F$-algebra 
generated by
the algebras $\F P$, $\F\W$ and $S(\h)$
with the following relations:
\begin{align*}
&s_i h =s_i (h) s_i-\bra\al_i\mid h\ket
\quad (i\in[1,n-1],\ h\in\h),\\
%\label{eq;trigrel2}
&s_i x^\eta s_i^{-1}= x^{s_i (\eta)}
\quad (i\in[1,n-1],\ \eta\in P),\\
%&s_i f s_i= s_i (f)
%\quad (i\in[1,n-1],\ f\in \F P),\\
%\label{eq;trigrel3}
&\left[ h, x^\eta\right]=\kappa\, \bra \eta\mid h\ket x^\eta
+\sum_{\al\in R^+}\bra \al\mid h\ket
{(x^\eta-x^{s_\al(\eta)})
\over 1-x^{-\al}} s_{\al}
\quad (h\in\h,\,\eta\in P).
\end{align*}
%where $e^\eta$ denote the element of $\F P\subset \F \affW$
%corresponding to $\eta\in P$.
\end{definition}

Define  {\it the degenerate affine Hecke algebra} of type $GL_n$
as the subalgebra of $\trigH$ generated by
$\F\W$ and $S(\h)$.
Observe that $\affH$ is a subalgebra of $\ratH$.

%Observe that the subalgebra $\F P\cdot \F\W$ is
%isomorphic to $\F \affW$.
%%
It is known due to Cherednik that
the natural multiplication map in $\trigH$ induces 
a linear isomorphism
$$\F P\* \F\W \* S(\h)\overset{\sim}{\rightarrow}\trigH.$$
%%%%%%%%%%%%%%%%%%%%%%%%%%%%%%%%%
Define the category 
$\O(\trigH)$ as the  full subcategory of $\trigHmod$ consisting of
the modules which are
%finitely generated over $\trigH$ and 
locally finite for $\UU$, i.e.,
a finitely generated $\trigH$-module $M$ is in $\O(\trigH)$
if and only if $\dim_\F \UU v<\infty$ for any $v\in M$.

%%%%%%%%%%%%%%%%%%%%%%%%%%%%%%%%%%%%%%%%%%%%%%%%%%%%%%%%%%%%%
%
%%%%%%%%%%%%%%%%%%%%%%%%%%%%%%%%%%%%%%%%%%%%%%%%%%%%%%%%%%%%%
Recall that 
the elements $\ech_i\in\trigH$ $(i\in[1,n])$
are realized as  Cherednik-Dunkl operators
(of trigonometric Dunkl operators) 
$$\kappa x_i\frac{\partial}{\partial x_i}+
\sum_{\al\in R^+}\brac{\al}{\ech_i}\frac{1}{1-x^{-\al}}(1-s_{\al})+i-1
$$
on the polynomial representation $\XX$
of $\trigH$ (\cite{Ch;unification}).
We call the subalgebra $S(\h)=\F[\ech_1,\dots,\ech_n]$
the Cherednik-Dunkl subalgebra.

For $\zeta\in\h^*$ and an $S(\h)$-module $M$,
define
\begin{align*}
M_\zeta&=\{v\in M \mid \ech_i v= \brac{\zeta}{\ech_i}v\ \forall i\in[1,n]\}
%\\ &=\{v\in M\mid fv=\chi_\zeta(f)v\ \forall f\in\UU\}.
\end{align*}
\def\wt{\h^*}
An element  $\zeta$ of $\wt$ is said to be a weight of $M$
if $M_\zeta\neq0$,
and an element of $M_\zeta$ a weight vector of weight $\zeta$.

Define $\Oss(\trigH)$
as the full subcategory of $\O(\trigH)$
consisting of those modules $M$ such that
$M=\oplus_{\zeta\in\wt}M_\zeta.$
%

%%%%%%%%%%%%%%%%%%%%%%%%%%%%%%%%%%%%%%%%%%%%%%%%%%%%%%%%%
\subsection{Tableaux representations}
%%%%%%%%%%%%%%%%%%%%%%%%%%%%%%%%%%%%%%%%%%%%%%%%%%%%%%%%%%
%%%%%%%%%%%%%%%%%%%%%%%%%%%%%%%%%%%%%%%%%%%%%%%%%%%%%%%%%%%%%
Now, we introduce representations of $\trigH$ associated with
 periodic diagrams.

Fix $\kappa\in\F$.
Let $m\in\Z_{\geq1}$ and put $\peri=(-m,\kappa-m)$.

For $\diag\in\dset_\peri^n$,
we denote by $\trigVdiag$ the space $\F\TabRC(\adiag)$
 of functions on the set $\TabRC(\adiag)$,
and denote by $\vec_T$ the element of $\trigVdiag$
corresponding to $T\in\TabRC(\adiag)$.
Namely, we set
$$\Vlm=\F\TabRC(\adiag)
=\bigoplus_{T\in\TabRC(\adiag)}\F \vec_T.$$
%%%%%%%%%%%%%%%%%%%%%%%%%%%%%%%%%%%%%%%%%%%%%%%%%%%%%%%%%%%%%%%%%%%
For $T\in\tab(\diag)$, define $\weight_{T}\in\h^*$ by
\begin{equation}
  \label{eq;weight_T}
\weight_T=\sum_{i\in[1,n]}\con_T(i)\e_i,
\end{equation}
where $\con_T$ is the content of $T$.
%%%%%%%%%%%%%%%%%%%%%%%%%%%%%%%%%%%%%%
\begin{lemma} \label{lem;irrlemma}
$\mathrm{(i)}$
Let $\diag\in\dset^n$ and $T\in\tabRC(\diag)$.
Let $M$ be an $\affH$-module and
suppose that $M_{\weight_T}\neq0$.
Then $\affH v \cap M_{\weight_S}\neq0$ for any $v\in M_{\weight_T}$ and
$S\in\tabRC(\diag)$.

\smallskip\noindent
$\mathrm{(ii)}$
Let $\diag\in\dset_\peri^n$ and $T\in\TabRC(\adiag)$.
Let $M$ be an $\trigH$-module and suppose that $M_{\weight_T}\neq0$.
Then $\trigH v\cap M_{\weight_S}\neq0$ for any $v\in M_{\weight_T}$ and $S\in\TabRC(\adiag)$.
\end{lemma}
\noindent{\it Proof.}
We will prove the statement (ii). 
The statement (i) follows similarly.

Let $v\in M_{\weight_T}$.
Let 
 $S\in\TabRC(\adiag)$ and take $w\in\affW$ such that $S=wT$.
We will prove that
$\trigH v\cap M_{\weight_S}=\trigH v\cap M_{\weight_{wT}}\neq 0$
by induction on the length $l(w)$ of $w$.
If $l(w)=0$ then $w=\pi^k$ and we have
$\pi^k v\in M_{\weight_{wT}}\setminus\{0\}$, and the statement holds.

Suppose that $l(w)=k$ and
the claim holds for all $z\in\affW$ 
with $l(z)<k$.
Take $i\in[0,n-1]$ and $z\in\affW$ such that $w=s_i z$ and 
$l(w)=l(z)+1$.
Then by Lemma~\ref{lem;weakorder}, we have $zT\in\TabRC(\adiag)$, and
by the induction hypothesis, we can find a non-zero element 
$v'$ in $\trigH v\cap M_{\weight_{zT}}$.
Put $c=\con_{zT}(i)-\con_{zT}(i+1)$. % and put $v'=(1-c_i s_i)v$.
Then is easy to check that  $(1+cs_i)v'\in M_{\weight_{wT}}$.
Note that $c\neq \pm 1$ by Lemma~\ref{lem;noncritical},
and hence $(1-c s_i)(1+cs_i)v'=(1-c^2)v'\neq0$.
Therefore $(1+cs_i)v'$ is a non-zero element in 
$\trigH v\cap M_{\weight_{wT}}$.
\qed
%%%%%%%%%%%%%%%%%%%%%%%%%%%%%%%%%%%%%%%%%%%%%%%
\begin{theorem}
$($cf. Theorem~3.16, Theorem 3.17 in \cite{SV}$)$
\label{th;tabrep_trig}
Let $\kappa\in\F$, $m\in\Z_{\geq1}$ and put $\peri=(-m,\kappa-m)$.
Let $\diag\in\dset_\peri^n$.

\smallskip\noindent
$\rm{(i)}$
There exists a unique $\trigH$-module
structure on $\Vlm$ such that
%\noindent
%${\rm(a)}$ The space  $\Vlm$ has a basis labeled by the set 
$\TabRC(\adiag)$:
\begin{align*}
\ech_i \vec_T&=\con_T(i) \vec_T\quad (i\in[1,n]),\\
\pi \vec_T&=\vec_{\pi T},\\
s_i \vec_T&=
\begin{cases}
\frac{1+\con_T(i)-\con_T(i+1)}{\con_T(i)-\con_T(i+1)}\vec_{s_i T}
-\frac{1}{\con_T(i)-\con_T(i+1)}\vec_T\
&\text{if }s_iT\in\TabRC(\adiag)\\
-\frac{1}{\con_T(i)-\con_T(i+1)}\vec_T\
&\text{if }s_iT\notin\TabRC(\adiag)
\end{cases}
\ \ (i\in[0,n-1]).
%\label{eq;Taction}
\end{align*}
%where
$($Note that $\con_T(i)-\con_T(i+1)\neq 0$
%which is non-zero
by Lemma~\ref{lem;noncritical}$)$.

\smallskip
\noindent
$\rm{(ii)}$
The $\Vlm$ admits a weight space decomposition with respect to
the %Cherednik-Dunkl 
subalgebra $S(\h):$
$$\Vlm
%\!\downarrow_{S(\h)}
=\bigoplus_{T\in\TabRC(\adiag)} \Vlm_{\weight_T},$$
and, moreover,  $\Vlm_{\weight_T}=\F \vec_T$ for all $T\in\TabRC(\adiag)$.

\smallskip
\noindent
$\rm{(iii)}$ The $\trigH$-module $\Vlm$ is irreducible.
\end{theorem}
\noindent{Proof.}
(i) The statement is proved by verifying the defining relation of $\trigH$
by direct calculations.

\noindent
(ii) Follows directly from Lemma \ref{lem;con_different}.
%For $S,T\in\TabRC(\adiag)$, it is easy to see that
%$S\neq T$ implies $\con_S\neq \con_T$.
%This implies the statement.

\noindent
(iii)
Take any $S,T\in\TabRC(\adiag)$
and take $w\in\affW$ such that $S=wT$.
By Lemma \ref{lem;irrlemma}, we have
$\vec_S\in \trigH \vec_T$.
This implies the irreducibility.
\qed

\medskip
%%%%%%%%%%%%%%%%%%%%%%%%%%%%%%%%%%%%%%%%%%%%%%%%%%%%%%%%%%%%%%%%%%
We call $\trigVdiag$ %(resp., $\affV(\diag)$)
the tableaux representation of $\trigH$ %(resp., $\affH$)
associated to $\adiag$. % (resp., $\diag$).

\begin{remark}
%Let $\kappa\in\Qgeq$.
%By modifying 
 For $\kappa\in\F\setminus\Q_{\leq0}$,
it can be shown that
 any irreducible module in $\Oss(\trigH)$
is isomorphic to a tableaux representation $\trigVdiag$
for some $m\in[1,n]$ and $\diag\in\dset_{\mkappa}^{*n}$.
A proof of this statement is given 
by modifying the proof of
\cite[Theorem 4.20]{SV}, where 
the corresponding theorem
 for the double affine Hecke algebra
is proved with the restriction 
$\kappa\in\Z_{\geq1}$ and $\diag\subset\Z\times\Z$.

The same statement has been given by
Cherednik in \cite{Ch;fourier, Ch;book}, where
 the condition for
two diagrams to give isomorphic representations
is also given.
\end{remark}
Similarly, 
putting $\affVlsm=\oplus_{T\in\tabRC(\diag)}\F \bar \vec_T$,
we have the following:
%%%%%%%%%%%%%%%%%%%%%%%%
\begin{proposition}[\cite{Ch;special_bases,Ram}]   
\label{pr;tabrep_aff}
 There exists a unique 
$\affH$-module structure on $\affV(\diag)$ such that
\begin{align*}
%\label{eq;Vaff1}
\ech_i \bar \vec_T&=\con_T(i) \bar \vec_T\quad (i\in[1,n]),\\
s_i \bar \vec_T&=
\begin{cases}
\frac{1+\con_T(i)-\con_T(i+1)}{\con_T(i)-\con_T(i+1)}\bar \vec_{s_i T}
-\frac{1}{\con_T(i)-\con_T(i+1)}\bar \vec_T\
&\text{if }s_iT\in\tabRC(\diag)\\
-\frac{1}{\con_T(i)-\con_T(i+1)}\bar \vec_T\
&\text{if }s_iT\notin\tabRC(\diag)
\end{cases}
\ \ (i\in[1,n-1]).
%\label{eq;Vaff2}
\end{align*}
Moreover, %$\dim \affV(\diag)_{\wtcon_T}=1$ for
%any $T\in\tabRC(\diag)$, and
$\affV(\diag)$ is an irreducible  $\affH$-module.
 \end{proposition}
%%%%%%%%%%%%%%%%%%%%%%%%%
%%%%%%%%%%%%%%%%%%%%%%%%%%%%%%%%%%%%%%%%%%%%%%%%%%%%%%%
\subsection{Restriction rule and plane partitions}
%%%%%%%%%%%%%%%%%%%%%%%%%%%%%%%%%%%%%%%%%%%%%%%%%%%%%%%
Let $\kappa\in\F$, $m\in\Z_{\geq1}$ and put $\peri=(-m,\kappa-m)$.
%Let $\peri\in\Z_{\leq-1}\times\F$ and 
Let $\diag\in\dset^n_\peri$.
We put
%\begin{alignat*}{2}
$\adiag^T=T^{-1}([1,n])\subset \adiag$
for $T\in\Tab_\peri(\adiag)$
as before.
Recall the surjection $\projmap:\TabRC(\adiag)\to\trigPPp(\diag)$
in Proposition \ref{pr;TabandPP},
for which we have $\adiag^T={\projmapop{T}}^{-1}(0)$.
%%%%%%%%%%%%%%%%%%%%%%%%%%%%%%%%%%%%%%%%%%%%%%
\begin{lemma}\label{lem;Vinc}
Let $T\in\TabRC(\adiag)$.

\noindent
$\mathrm{(i)}$
Any $S\in\tabRC(\adiag^T)$ can be extended to
a standard tableau $\hat{S}$ on $\adiag$ by 
setting $\hat{S}(u+j\peri)=S(u)-jn$ $(u\in\adiag^T,\ j\in\Z)$.

\noindent
$\mathrm{(ii)}$
The assignment $\bar \vec_S\mapsto \vec_{\hat{S}}$ 
$(S\in\tabRC(\adiag^T))$ gives an embedding 
 $\affV(\adiag^T)\to\trigVdiag$
as an $\affH$-module.
\end{lemma}
%%%%%%%%%%%%%%%%%%%%%%%%%%%%%%%%%%%%%%%%%%%%%%%%%%
\noindent{\it Proof.}
(i) Let $S\in\tabRC(\adiag^T)$.
It is obvious that $S$ satisfies the condition (T1) in Definition \ref{def;tableaux}.

We shall check the condition (T3).
Let
$(a,b)\in\adiag^T\subset \adiag$.
It is enough to show that $S(a,b)<S(a+k+1,b+k)$
for any $(a+k+1,b+k)\in\adiag\setminus\adiag^T$
with $k\in\Z_{\geq0}$.

We have
$\projmapop{T}(a,b)=0$ as $\adiag^T=\projmapop{T}^{-1}(0)$, and we have
$\projmapop{T}(a+k+1,b+k)\geq1$ as $\projmapop{T}\in\trigPP_\peri(\diag)$
by Proposition \ref{pr;TabandPP}.
Note that there exists $w\in\W$ such that $S=wT$ on
$\adiag^T$.
Hence we have $\projmapop{S}=\projmapop{T}$.
This implies
$\hat{S}(a,b)<\hat{S}(a+k+1,b+k)$.

The condition (T2) can be checked similarly.
Therefore $\hat{S}\in\TabRC(\adiag)$.

\noindent
(ii) It follows from the first statement that
the defining relations 
of $\affV(\adiag^T)$ in Proposition \ref{pr;tabrep_aff}
hold for
$\{\vec_{\hat{S}}\}_{S\in\tabRC(\adiag^T)}$.
This implies $\oplus_{S\in\tabRC(\adiag^T)}\F \vec_{\hat S}
\cong\affV(\adiag^T)$,
and the statement follows.
\qed

\medskip
We have the following simple decomposition formula
for the restriction $\trigVdiag\!\downarrow_{\affH}$ 
as an $\affH$-module (cf. \cite[Theorem 3.7.3]{Ch;book}).
%
%%%%%%%%%%%%%%%%%%%%%%%%%%%%%%%%%%%%%%%%%%%%%%%%%%%%%%%%
\begin{theorem}\label{th;restriction_trig}
%$($cf. \cite[Theorem 3.7.3]{Ch;book}$)$
Let $\kappa\in\F$ and $m\in\Z_{\geq1}$.
Let  $\diag\in\dset^n_\peri$ with $\peri=\mkappa$. 
Then
  $$\trigVdiag\!\downarrow_{\affH}
=\bigoplus_{\p\in\trigPPp(\diag)}
\affV(\p^{-1}(0)).$$
\end{theorem}
%%%%%%%%%%%%%%%%%%%%%%%%%%%%%%%%%%%%%%%%%%%%%%%%%%%%%%%%%%
\noindent{\it Proof.}
Let $T\in \TabRC(\adiag)$.
By Lemma \ref{lem;Vinc},
we have $\trigVdiag\!\downarrow_{\affH}
=\sum_{T\in\TabRC(\adiag)}\affV(\adiag^T)$.

Let $T,S \in\TabRC(\adiag)$, then 
$\affV(\adiag^T)\cap\affV(\adiag^S)=\{0\}$
or $\affV(\adiag^T)=\affV(\adiag^S)$ %as subspaces of $\trigVdiag$,
because they %$\affV(\adiag^T)$ 
are irreducible. 
Moreover, it can be seen easily that
\begin{align*}
\affV(\adiag^T)=\affV(\adiag^S)
&\Leftrightarrow \vec_T\in\affV(\adiag^S)=\F\W\cdot \vec_S\\
&\Leftrightarrow T=wS\hbox{ for some }w\in\W.
\end{align*}
Now, the statement follows from 
$\adiag^T={\projmapop{T}}^{-1}(0)$ and Proposition \ref{pr;TabandPP}
\qed
%%%%%%%%%%%%%%%%%%%%%%%%%%%%%%%%%%%%%%%%%%%%%%%%%%
\subsection{$\W$-invariant subspace}
%%%%%%%%%%%%%%%%%%%%%%%%%%%%%%%%%%%%%%%%%%%%%%%%%%
%\medskip
Define the elements $\unip_\pm$ in $\F\W$ by
\begin{align}
&\unip_+=\frac{1}{n!}\sum_{w\in\W}w,\\  
&\unip_-=\frac{1}{n!}\sum_{w\in\W}(-1)^{l(w)}\,w.
\end{align}
The space $\unip_+ M$ (resp., $\unip_- M$)
is the $\W$-invariant (resp.,  $\W$-anti-invariant) subspace of $M$:
$\unip_{\pm}M=\{v\in M\mid wv=(\pm 1)^{l(w)} v\ \forall w\in\W \}$.
The algebra $S(\h)^\W$ acts on $\unip_{\pm}M$.
We will give a decomposition of $\unip_{\pm}\trigVdiag$ as
an $S(\h)^\W$-module.

%%%%%%%%%%%%%%%%%%%%%%%%%%%%%%%%%%%%%%%%%%%%%%%%%%%%%
\begin{definition}
A skew diagram $\diag$
is said to be {\it linked}
 if 
there exist $(a,b)\in\diag$ and  $(a+k+1,b+k)\in\diag$
with $k\in\Z_{\geq0}$.

A skew diagram $\diag$
is said to be {\it unlinked} if it is not linked.
\end{definition}
%%%%%%%%%%%%%%%%%%%%%%%%%%%%%%%%%%%%%%%%%%%%%%%%%%%%%%%
\begin{remark}
When $\diag\subset\Z\times\Z$, an unlinked diagram $\diag$ 
is also called a horizontal strip.
\end{remark}
%%%%%%%%%%%%%%%%%%%%%%%%%%%%%%%%%%%%%%%%%%%%%%%%%%%%%%%%%%%%%%%%
%%%%%%%%%%
A proof of the following proposition is given in the next section.
%%%%%%%%%%%%%%%%%%%%%%%%%%%%%%%%%%%%%%%%%%%
\begin{proposition} \label{pr;unlinked}
Let $\diag$ be a skew diagram. Then
\begin{align*}
&\dim_\F\,  \unip_+\! \affV(\diag)=
\begin{cases}&1\ \ \hbox{ if } \diag\hbox{ is unlinked},\\
             &0\ \ \hbox{ if } \diag\hbox{ is linked},
\end{cases}\\
&\dim_\F \, \unip_-\! \affV(\diag)=
\begin{cases}&1\ \ \hbox{ if } \diag'\hbox{ is unlinked},\\
             &0\ \ \hbox{ if } \diag'\hbox{ is linked},
\end{cases}
\end{align*}
where $\diag'$ is the conjugate of $\diag:$
$\ \diag'=\{(a,b)\in\Z\times\F \mid (b,a)\in\diag\}$.
\end{proposition}
%%%%%%%%%%%%%%%%%%%%%%%%%%%%%%%%%%%%%%%%%%%%%%%%%%%%%%%%%%%

%\medskip

For $\zeta\in\h^*$, 
 let 
 $\chi_\zeta:S(\h)^\W \to\F$
denote the character corresponding to the image $\W\zeta$ of $\zeta$ 
under the natural projection
$\h^*\to\W\backslash \h^*$.

Let $\peri=(-m,\kappa-m)$ with $m\in\Z_{\geq1}$.
Recall the isomorphism 
$\projmap:\Worb{\TabRC}(\adiag)\isomto\trigPPp(\diag)$.
For $\p\in\trigPPp(\diag)$, define $\chi_\p=\chi_{\weight_T}$, where
$T$ is any standard tableau such that $\projmapop{T}=\p$.

\begin{lemma}
  Let $\p,\mathfrak{q}\in\trigPPp(\diag)$.
%$T,S\in\TabRC(\adiag)$.
Then $\chi_\p=\chi_{\mathfrak{q}}$ if and only if
$\p=\mathfrak{q}$.
\end{lemma}
\noindent{\it Proof.}
We have the map
$\TabRC(\adiag)\to\h^*$
given by $T\mapsto \weight_T$.
This map is injective by Lemma \ref{lem;con_different}.
Moreover, it follows from Lemma \ref{lem;con_weq}
that $w(\weight_T)=\weight_{wT}$.
%, namely this map is
%$\W$-equivaliant.
Hence it factors
the injection
$\Worb{\TabRC}(\adiag)\to \W\backslash\h^*$.
Therefore we have an injection 
$\trigPPp(\diag)\to \W\backslash\h^*$.
%It is obvious from the definition that $\chi_\p=\chi_\zeta$ 
%if $\tilde\projmapop{\p}=\W\zeta$.
This implies the statement.\qed

\medskip
For a character $\chi$ of $S(\h)^\W$, put  
$$\unip_\pm\trigVdiag^{\chi}=\{ v\in\unip_\pm \trigVdiag\mid
h v=\chi(h) v\ \ \forall h\in S(\h)^\W \}.$$
For $T\in\TabRC(\adiag)$,  observe that
the vectors $\unip_+ \vec_T\in\trigVdiag$ and   $\unip_- \vec_T\in\trigVdiag$ 
have
the character $\chi_{\projmapop{T}}=\chi_{\weight_T}$.
%%%
%%%%%%%%%%%%%%%%%%%%%%%%%%%%%%%%%%%%%%%%%%%%%%%%%%%%%%%%%%
\begin{theorem}\label{th;VandPP}
Let $\kappa\in\F$ and $m\in\Z_{\geq1}$, and
let $\diag\in\dset_\peri^n$ with $\peri=\mkappa$. Then

$$\unip_+\trigVdiag=\bigoplus_{\p\in\trigPPCp(\diag)}
\unip_+\trigVdiag^{\chi_\p},$$
$$\unip_-\trigVdiag=\bigoplus_{\p\in\trigPPRp(\diag)}
\unip_-\trigVdiag^{\chi_\p}.$$
Moreover, 
 $\dim_\F\; \unip_+\trigVdiag^{\chi_\p}=1$
for all $\p\in\trigPPCp(\diag)$,
and  
$\dim_\F\; \unip_-\trigVdiag^{\chi_\p}=1$
for all $\p\in\trigPPRp(\diag)$.
\end{theorem}
\noindent{\it Proof.}
By Theorem \ref{th;restriction_trig} and Proposition \ref{pr;unlinked},
we have
$$\unip_+\trigVdiag
=\bigoplus_{\p\in\trigPPp(\diag),\ \p^{-1}(0):unlinked}
\unip_+ \affV(\p^{-1}(0))$$
To prove the formula for $\unip_+\trigVdiag$, it is
enough to show that
 $\p^{-1}(0)$ is unlinked if and only if $\p\in\trigPPCp(\diag)$.
The ``if'' part is obvious. 
Let us prove the opposite implication.
Suppose that $\p$ is not column strict.
Then there exists $u_1=(a,b)\in\adiag$ and $u_2=(a+k+1,b+k)\in\adiag$ 
with $k\in\Z_{\geq0}$
such that $\p(u_1)=\p(u_2)$. %Denote this number by $p$.
Then $\p(u_1+\p(u_1)\peri)=\p(u_2+\p(u_1)\peri)=0$.
Hence $u_1+\p(u_1)\peri$ and $u_2+\p(u_1)\peri$ belongs 
to $\p^{-1}(0)$,
and this implies that $\p^{-1}(0)$ is linked.

Therefore we have proved the formula for $\unip_+\trigVdiag$.
The formula for $\unip_+\trigVdiag$ follows similarly.\qed
%%%%%%%%%%%%%%%%%%%%%%%%%%%%%%%%%%%%%%%%%%%%%%%%%%%%%%%%%%%%%%%%%%%%%%%%

\medskip
%By Proposition \ref{pr;tabres_PPC}, we have the following 
%decomposition formula, which
The following decomposition formula
was conjectured in \cite[Conjecture 6.2.6]{AST}.
%
%%%%%%%%%%%%%%%%%%%%%%%%%%%%%%%%%%%%%%%%%%%%%%%%%%%%%%%%%%%%%%%%%%
\begin{theorem}\label{th;sp_decomp}
Let  $\kappa\in \Z_{\geq1}$, $m\in\Z_{\geq1}$ and
$\diag\in\dsetint^{n}_\peri$
with $\peri={\mkappa}$.
Then 
  $$\unip_+\trigVdiag=
\bigoplus_{\zeta\in P^-, \ T\in\restab_{\peri}(\diag)}
\unip_+\trigVdiag^{\chi_{\zeta+\ccwt_T}},$$
Moreover 
$\unip_+\trigVdiag^{\chi_{\zeta+\ccwt_T}}
=\F \unip_+ \vec_{t_{\zeta+\ccwt_T} T}$ and it is non-zero
 for all $\zeta\in P^-$
and $T\in\restab_\pri(\diag)$.
\end{theorem}
\noindent
{\it Proof.}
Follows from Proposition \ref{pr;tabres_PPC}
and Theorem \ref{th;VandPP}
noting that 
$\chi_{\projmapop{t_{\zeta+\ccwt_T} T}}=\chi_{\zeta+\ccwt_T}$
for $T\in\tabRC_\peri(\diag)$.\qed
%%%%%%%%%%%%%%%%%%%%%%%%%%%%%%%%%%%%%%%%%%%%%%%%%%%%%%%%%
\subsection{Proof of Proposition \ref{pr;unlinked}}
%%%%%%%% Parabolic Subgroup%%%%%%%%%%%%%%%%%%%%%%%%%%%%%%
%
First, we prepare notations concerning cosets of the Weyl group.

For $w\in\W $, set 
$ R(w)= R^+\cap w^{-1}  R^-,$
where $ R^-=  R\setminus   R^+$.

Let $\nu=(\nu_1,\dots,\nu_m)\models n$ (a composition of $n$).
% be a composition of $n$,i.e., $\sum_{i\in[1,m]}\nu_i=n$. 
Put $I(\nu)=[1,n]
\setminus \{\nu_1,\nu_1+\nu_2,\dots,\sum_{k\in[1,m]}\nu_k\}$
and define 
\begin{align*}
%\Pi_\nu&=\{\al_i\mid i\in I(\nu)\}\subseteq\Pi,\\
\W_\nu&=\bra s_i\mid i\in I(\nu)\ket\subseteq W.
%\label{eq;parabolicsub}
%\\
%R_\nu^+&=\{\al\in R^+\mid s_\al\in W_\nu\}.
\end{align*}
The subgroup $\W_\nu$ is called the parabolic subgroup 
(or Young subgroup)
associated with $\nu$.
%corresponding %to $\Pi_\nu$.
%
Define 
\begin{align}\label{eq;coset_rep_W}
& \W^{\nu}=\left\{w\in \W\mid R(w)\cap R_\nu^+=\emptyset
\right\},
%\label{eq;coset_rep_affW}
%& \affW^{\nu}=\left\{w\in \affW\mid R(w)\cap  R_\nu^+=\emptyset
%\right\}.
\end{align}
where $R_\nu^+=\{\al\in R^+\mid s_\al\in W_\nu\}.$
As is well-known, 
 $\W^{\nu}$ %(resp., $\affW^{\nu}$) gives a 
complete set of representatives of the coset $\W/\W_\nu$.
%(resp., $\affW/\affW_\nu$).

%%%%%%%%%%%%%%%%%%%%%%%%%%%%%%%%%%%%%%%%%%%%%%%%%%%%%%%%%
\medskip
\noindent{\it Proof of Proposition \ref{pr;unlinked}.}
%%%%%%%%%%%%%%%%%%%%%%%%%%%%%%%%%%%%%%%%%%%%%%%%%%%%%%%%%%%%
Suppose that $\diag$ is linked.
Then there exist
$(a,b)\in\diag$ and 
$(a+k+1,b+k)\in\diag$ with $k\in\Z_{\geq0}$.
We take $k$ as small as possible.
Then it is easy to see that
there exists 
$T\in\tabRC(\diag)$ such that 
$T(a+k+1,b+k)=T(a,b)+1$.

Put $i=T(a,b)$.
%For such $T\in\tabRC(\diag)$, we have
Then $s_iT\notin\tabRC(\diag)$,
and
%such that $T(a+k+1,b+k)=T(a,b)+1$, 
it follows from the definition of the
action of $s_i$ (Theorem \ref{th;tabrep_trig})
that $s_i \vec_T=-\vec_T$,
and hence
$$\unip_+ \vec_T=\frac{1}{n!}\sum_{w\in\W} w \vec_T=
\frac{1}{n!}\sum_{w\in\W^{\nu(i)}} w(1+s_i) \vec_T=0.$$
where $\nu(i)$ denotes the composition $(i,n-i)$,
and $\W^{\nu(i)}$ denotes the set
of coset representatives of $\W/\W_{\nu(i)}$.
%(see \eqref{eq;coset_rep_W}).

Since $ \affV(\diag)= \affH \cdot \vec_T=\F\W \cdot \vec_T$,
we have
$\unip_+ \affV(\diag)=
\unip_+\, \F\W \cdot \vec_T=\F \unip_+ \vec_T
=0.$

Suppose that $\diag$ is unlinked
and take $\lm,\mu\in\F^m$ such that $\diag=\lm/\mu$. 
Let $M$ denote the $\F\W$-module
$\F\W\otimes_{\F
\Weyl_{\lm-\mu}}\F\one_{\lm-\mu}$,
where  $\F \one_{\lm-\mu}$ 
denotes the trivial module of the parabolic subgroup
$\Weyl_{\lm-\mu}$.

Observe that $s_i\bar \vec_{\rbtab}=\bar \vec_{\rbtab}$
for any $s_i\in\Weyl_{\lm-\mu}$.
Here, $\rbtab$ is the row reading tableau on
 $\diag$ in \eqref{eq;rbtab}.
Hence there exists an $\F\W$-homomorphism 
$$M
{\to}\affV(\diag)$$
such that $\one_{\lm-\mu}\mapsto \bar{v}_\rbtab$.
Note that we have 
$\affV(\diag)=\affH \cdot \bar \vec_{\rbtab}
=\F\W \cdot\bar \vec_{\rbtab}$ as
$\bar \vec_{\rbtab}$ is a weight vector with respect to $S(\h)$.
Hence the map $M\to\affV(\diag)$ is surjective.

It is easy to see that 
$$\sharp \tabRC(\diag)=\frac{n!}{(\lm_1-\mu_1)!\dots(\lm_m-\mu_m)!}=
\sharp \W/\Weyl_{\lm-\mu}$$ 
for an unlinked diagram $\diag$.
Therefore $\dim \affV(\diag)=\dim M$,
which implies that the map above
is an $\F\W$-isomorphism.
Therefore $\unip_+ \affV(\diag)=\unip_+ M=\F\unip_+\one_{\lm-\mu}$
and it is one-dimensional.
This completes the proof of the statement for $\unip_+\affV(\diag)$. 
The statement for $\unip_-\affV(\diag)$ is proved similarly.
\qed
%%%%%%%%%%%%%%%%%%%%%%%%%%%%%%%%%%%%%%%%%%%%%%%%%%%%%%%%%%%%%%%%%%%%%%
%

%%%%%%%%%%%%%%%%%%%%%%%%%%%%%%%%%%%%%%%%%%%%%%%%%%
\section{Application to the rational Cherednik algebra}
%%%%%%%%%%%%%%%%%%%%%%%%%%%%%%%%%%%%%%%%%%%%%%%%%%%%

We will apply 
the combinatorial method developed in the previous
section to the study of representations
of the rational Cherednik algebra $\ratH$.

Representations of $\trigH$ with weight decomposition with 
respect to the Cherednik-Dunkl subalgebra
are constructed as tableaux representations
by giving the action of generators of $\trigH$ on
a basis vector labeled by standard tableaux.
On the other hand, such construction has not been known for $\ratH$.
But it will be shown that any irreducible $\ratH$-modules of 
the corresponding class can be realized as a subspace 
of tableaux representations of $\trigH$ using
the functor given in \cite{Su;ratandtrig},
which relates the representation theory of
$\ratH$ and $\trigH$.

%%%%%%%%%%%%%%%%%%%%%%%%%%%%%%%%%%%%%%%%%%%%%%%%%%%%%%%%
\subsection{Rational Cherednik algebra}
%%%%%%%%%%%%%%%%%%%%%%%%%%%%%%%%%%%%%%%%
We put
$\Preg=\oplus_{i\in[1,n]}\Z_{\geq 0}\e_i$
as before, which is 
the semi subgroup of $P$ generated by $\{\e_1,\dots,\e_n\}$.
Recall also that $\F \Preg=\F[\underline x]$ under our identification.
We define
$$\affWreg=\Preg\rtimes\W\subset \affW,$$
%and recall that $\affWreg$
which is a semigroup (with unit) with 
generators $\{\pi,s_1,s_2.\dots,s_{n-1}\}$.
\begin{definition}
For $\kappa\in\F$,   
the {\it the rational Cherednik algebra} $\ratH$ of type $GL_n$
is the unital associative $\F$-algebra 
generated by
$\F\affWreg$ and $\YY =\F[y_1,\dots,y_n]$
with the following relations:
\begin{align*}
&s_i y_j =y_{s_i(j)} s_i
\quad (i\in[1,n-1],\ j\in[1,n]),\\
&\left[ y_i, x_i\right]=
\begin{cases}
\kappa+\sum_{k\neq i} s_{ik}\ \ &(i=j)\\
-s_{ij}\ \ &(i\neq j)
\end{cases}
\quad (i,j\in[1,n]).
\end{align*}
\end{definition}

The natural map  gives a linear isomorphism
$$\F\Preg\* \F\W \* \YY\overset{\sim}{\rightarrow}
\ratH.$$
%%%%%%%%
Put $\Xch_+=\sum_{i\in[1,n]}\xch_i\Xch\subset \ratH$.
An element $v$ of an $\ratH$-module
is said to be $\Xch_+$-nilpotent
if there exists $k\in\Z_{>0}$ such that
$(\xch_i)^k v=0$ for all $i\in[1,n]$.
An $\ratH$-module $M$ is said to be
locally nilpotent for $\Xch_+$, if
any element of $M$ is  $\Xch_+$-nilpotent.

Define $\O(\ratH)$ as the full subcategory of $\ratHmod$ consisting of
those modules which are locally nilpotent for $\Xch_+$.
The following statement is easy to show but will play an important role.
%%%%%%%%%%%%%%%%%%%%%%%%%%%%%%%%%%%%%%%%%%
\def\emb{\iota}
\def\embinv{\jmath}
\begin{proposition}[\cite{Su;ratandtrig}]
\label{pr;embedding}
%${\rm (i)}$ 
There exists an algebra embedding
$\emb : \ratH \to \trigH$
such that
\begin{align*}
\emb(w)&=w \quad (w\in\W),\ \
\emb(x_i)=\e_i \quad (i\in[1,n]),\\
\emb(\xch_i)&=x_i^{-1}\left(\ech_i-\sum\nolimits_{1\leq
j<i}s_{ji}\right)\quad
(i\in[1,n]).
\end{align*}
%\smallskip
%\noindent
%${\rm (ii)}$
Moreover, the embedding $\emb:\ratH\to\trigH$
is extended to the algebra isomorphism
 $\XL\*_{\X} \ratH\isomto
\trigH$.
\end{proposition}
%%%%%%%%%%%%%%%%%%%%%%%%%%%%%%%%%%%%%%%%%%%%%%%%%%%%%%%
In the sequel, 
we often identify $\ratH$ with the subalgebra
of $\trigH$ generated by
$\F\affWreg$ and $y_i=x_i^{-1}(\ech_i-\sum\nolimits_{1\leq
j<i}s_{ji})$ $(i\in[1,n])$
via the embedding $\emb$.

Under this identification, $\ech_i\in\trigH$ is also
contained in $\ratH$, and it is expressed
as  $\ech_i=x_iy_i+\sum_{j<i}s_{ji}$ 
in terms of the generators of
$\ratH$.
%We call $\ech_i$ Cherednik-Dunkl elements.

Define $\Oss(\ratH)$
as a full subcategory of $\O(\ratH)$ 
consisting of those modules $M$ such that
$M=\oplus_{\zeta\in\wt}M_\zeta,$
where $M_\zeta=\{v\in M\mid \ech_i v=\brac\zeta{\ech_i} v \ \forall i\in[1,n]\}$
as before.

%%%%%%%%%%%
The following lemma is easily checked.
\begin{lemma}\label{lem;invo}
$\mathrm{(i)}$ 
There exists an algebra isomorphism $\sigma:\trigH\to \trigHm$
such that
\begin{align*}
&\sigma(x_i)=(-1)^{n-1}x_i,\ \ 
\sigma(\ech_i)=-\ech_i\ \ (i\in[1,n]),\\
&\sigma(s_i)=-s_i\ \ (i\in [0,n-1]),\
\end{align*}
and its restriction to $\ratH$ gives an algebra isomorphism
$\sigma: \ratH\isomto \H_{-\kappa}$.
%such that  $\sigma(\xch_i)=(-1)^n\xch_i$ $(i\in[1,n])$. %%%removed

\smallskip\noindent
$\mathrm{(ii)}$ 
The algebra isomorphism $\sigma$ 
%in Lemma \ref{lem;invo}
induces categorical equivalences
\begin{align*}
 & \O(\trigH)\isomto\O(\trigHm),\ \ 
\Oss(\trigH)\isomto\Oss(\trigHm),\\ 
& \O(\ratH)\isomto\O(\H_{-\kappa}),\ \ 
\Oss(\ratH)\isomto\Oss(\H_{-\kappa}).
\end{align*}
\end{lemma}

Consider the induction functor $\ratHmod\to\trigHmod$ given by 
$$M\to\trigH\otimes_{\ratH}M.$$
%from the category $\ratHmod$ to the category $\trigHmod$.
We denote its restriction to
$\O(\ratH)$ by $\Ind$.
Then it turns out that
 $\Ind$ gives an exact functor into $\O(\trigH)$
(\cite[Corollary 3.4, Proposition 4.2]{Su;ratandtrig}).
%%%%%%%%%%%%%%%%%%%%%%%%%%%%%%%%%%%%%%%%%%%%%%%%%%%%%%%%

For an $\trigH$-module $N$,
denote by $\nil(N)$ 
the subspace of $N$ consisting of
the $\Xch_+$-nilpotent elements:
\begin{equation}\label{eq;nil}
\nil(N)=\{v\in N\mid
v\hbox{ is }\Xch_+\hbox{-nilpotent}
\}.
\end{equation}
It follows that $\nil(N)$ is an $\ratH$-submodule
of $ N\!\!\downarrow_{\ratH}$ and it is finitely generated over $\ratH$
(\cite[Lemma 4.4]{Su;ratandtrig}).
%\end{lemma}
Therefore we have the functor
$\nil:\O(\trigH)\to\O(\ratH)$.
%%%%%%%%%%%%%%%%%%%%%%%%%%%%%%%%%%%%%%%%%%%%%%%%%%%%%%%%%%%%%%%%%%%%
\begin{theorem}$($\cite[Section 
6]{Su;ratandtrig}$)$\label{th;fullyfaithful}
$\mathrm{(i)}$ The functor $\nil$ is  the right adjoint functor of %the 
functor 
$\Ind$.
Moreover, $\nil\circ\Ind(M)\cong M$ for any $M\in\O(\ratH)$.
%\end{proposition}
%
%\begin{theorem}[\cite{S;ratandtrig}]\label{th;fullyfaithful}

\smallskip\noindent
$\mathrm{(ii)}$  The functor
$\Ind:\O(\ratH)\to\O(\trigH)$
is exact and fully-faithful.
\end{theorem}

%%%%%%%%%%%%%%%%%%%%%%%%%%%%%%%%%%%%%%%%%%%%%%%%%%%%%%%%%
\subsection{The category $\Oss$ for the rational algebra}
%%%%%%%%%%%%%%%%%%%%%%%%%%%%%%%%%%%%%%%%%%%%%%%%%%%%%%%%%
%

Let $m\in\Z_{\geq1}$ and $\kappa\in\Qgeq$,
and put $\peri=\mkappa$.

Let $\part(m,n)$ denote the set of partitions of $n$ consisting
$m$ nonzero components:
\begin{equation}
\part(\m,n)=
\{\lm=(\lm_1,\dots,\lm_m)\in(\Z_{\geq1})^m \mid\lm\vdash n\}
\end{equation}
Set 
$\part_\kappa(\m,n)=\part(\m,n)$ for $\kappa\in\F\setminus\Q$, and set
\begin{equation}
\part_\kappa(\m,n)=
\{\lm\in \part(\m,n)\mid
s-\m-\lm_1+\lm_\m\in\Z_{\geq0}\}
\end{equation}
for 
$\kappa=s/r\in\Q_{>0} \hbox{ with } s,r\in\Z_{>0},\ (s,r)=1.$

%Observe that $\part(\m,n)\cap \Domk =\part_\kappa(\m,n)$.

We identify a partition $\lm\in \part(\m,n)$ with
the associated diagram
$\{(a,b)\in\Z\times\Z\mid a\in[1,m], b\in[1,\lm_a]\}$.
By this identification,
$\part(m,n)$ and $\part_\kappa(m,n)$ are 
 thought as a subset of
$\dset^n$ and $\dset_\peri^n$ respectively.

\def\trigValm{\trigV(\alm)}

Recall the functor 
$\nil:\O(\trigH)\to\O(\ratH)$,
and consider its image of the tableau representation
$\trigValm$ for $\lm\in\part_\kappa(m,n)$.
By Theorem \ref{th;fullyfaithful},
%and Proposition \ref{pr;OssOss},
$\nil(\trigValm)$ is irreducible or zero.  Moreover, 
%$\nil(\trigValm)$ 
it admits a weight decomposition with respect 
to the Cherednik-Dunkl subalgebra $S(\h)$,
 since it is an $S(\h)$-submodule of $\trigValm$. Hence 
$\nil(\trigValm)$
belongs to $\Oss(\ratH)$.
%%%%%%%%%%%%%%%%%%%%%%%%%%%%%%%%%%%%%%%%%

\begin{remark}
$($See \cite[Proposition 8.1]{Su;ratandtrig}.$)$
%\label{pr;OssOss}
More generally, it holds for any $M\in\O(\ratH)$
that  
$$M\in\Oss(\ratH)\Leftrightarrow\Ind(M)\in \Oss(\trigH).$$
\end{remark}

%%%%%%%%%%%%%%%%%%%%%%%%%%%%%%%%%%%%%%%%%
%%%%%%%%%%%%%%%%%%%%%%%%%%%%%%%%%%%%%%%%%%%%%%%%%%%%%%%%%%%%%%%%%%%%%

%The following result is obtained
\def\sp{S}
For a partition $\lm\vdash n$,
let $\sp_\lm$ denote the corresponding irreducible $\W$-module.
%It is known
%an $\F\W$-module $E$, 
We define the {\it standard module} of $\ratH$
associated with $\sp_\lm$ by
$$\ratst(\lm)=\ratH\otimes_{\F\W\cdot \YY}\sp_\lm,$$
where we let $\YY$ act on $\sp_\lm$ through the augmentation map
$\YY\to\F$ given by $y_i\mapsto 0$  $(i\in[1,n])$.

Let $\Irr\O(\ratH)$ (resp., $\Irr\Oss(\ratH)$) denote
the set of isomorphism classes
of the irreducible modules in $\O(\ratH)$
(resp., $\Oss(\ratH)$).
%Recall that $S(\h)\subset\affH\subset\ratH\subset\trigH$.
%
%%%%%%%%%%%%%%%%%%%%%%%%%%%%%%%%%%%%%%%%%%%%%%%%%%%%%%%%%%%%%%%%%%%%%%%%%%%%%

\begin{proposition}
  $($\cite{DO,GGOR} %\cite[Proposition 2.11]{GGOR}
$)$\label{pr;GGORclassification}
Let $\kappa\in\F\setminus\{0\}$.

For $\lm\in\sqcup_{m\in[1,n]}\part(m,n)$, 
the standard module $\ratst(\lm)$ has a unique simple quotient module,
which is denoted by $\ratL(\lm)$.
Moreover,
the assignment $\lm\mapsto \ratL(\lm)$ gives
a one-to-one correspondence
$$\bigsqcup_{\m\in[1,n]}\part(\m,n)\isomto\Irr\O(\ratH).$$
\end{proposition}
%%%%%%%%%%%%%%%%%%%%%%%%%%%%%%%%%%%%%%%%%%%%%%%%

%%%%%%%%%%%%%%%%%%%%%%%%%%%%%%%%%%%%%%%%%%%%%%
\begin{proposition}\label{pr;L=nil}
%$\mathrm{(i)}$
Let $\kappa\in\F\setminus\Q_{\leq0}$ 
and $m\in\Z_{\geq1}$. 
For $\lm\in\part_\kappa(m,n)$, it holds that
$\ratL(\lm)\cong\nil(\trigValm)$, where  $\peri=(-m,\kappa-m)$.
% isthe unique simple quotient of $\ratst(\lm)$.
\end{proposition}
\noindent
{\it Proof.}
Recall that  
$\affV(\lm)$ is embedded into $\trigValm$
via $\bar \vec_T\mapsto \vec_T$ as an $\affH$-submodule.
Obviously, (the image of) $\affV(\lm)$
is locally $\YY_+$-nilpotent, and hence
$\affV(\lm)\subset \nil(\trigValm)$.
Remark that 
for $\lm\in\part(m,n)$, 
the restriction $\affV(\lm)\!\!\downarrow_{\F\W}$ is 
irreducible as an $\F\W$-module, and it is isomorphic to $\sp_\lm$.
Therefore we have a surjective homomorphism
$\ratst(\lm)\to \nil(\trigValm).$
This implies that $\nil(\trigValm)$ 
is the simple quotient of $\ratst(\lm)$.
\qed

\medskip
We have the following result for $\Oss(\ratH)$,
whose proof will be given in Appendix.
%%%%%%%%%%%%%%%%%%%%%%%%%%%%%%
\begin{theorem}$($cf. \cite[Theorem 8.2]{Su;ratandtrig}$)$
\label{th;classification_rat}
  Let $\kappa\in\C\setminus\Q_{\leq 0}$.
%\smallskip\noindent
%$\mathrm{(ii)}$
The assignment $\lm\mapsto \ratL(\lm)$ gives
a one-to-one correspondence
$$\bigsqcup_{\m\in[1,n]}\part_\kappa(\m,n)\isomto\Irr\Oss(\ratH).$$
\end{theorem}
%%%%%%%%%%%%%%%%%%%%%%%%%%%%%%%%%%%%%%%%%%%%%%%%
\begin{remark}
(i) The classification result in Theorem \ref{th;classification_rat}  
for $\kappa\in \Q_{<0}$
can be obtained via the isomorphism $\ratH\cong{\mathcal H}_{-\kappa}$.

\smallskip\noindent 
(ii) In the case where $\kappa\in\Z\setminus\{0\}$, the classification 
%in
%Theorem \ref{th;classification_rat}
 is given in \cite{Su;ratandtrig}
as a consequence of the classification theorem for $\Oss(\trigH)$
(\cite[Theorem 6.5]{Ch;fourier}\cite[Corollary 4.23]{SV}).
The statement for general $\kappa$ is also stated in 
\cite[Remark 8.3]{Su;ratandtrig}
without a precise proof.
It should be also mentioned that, in the proof given in Appendix,  
the classification result for $\trigH$, which is more complicated, is not 
used.  
\end{remark}
%%%%%%%%%%%%%%%%%%%%%%%%%%%%%%%%%%%%%%%%%%%%%%%%%%%%%%%
\subsection{Tableaux description for irreducible $\ratH$-modules}
%%%%%%%%%%%%%%%%%%%%%%%%%%%%%%%%%%%%%%%%%%%%%%%%%%%%%

We will give a weight decomposition 
of $\ratVlm=\nil(\trigValm)$ with respect to 
the subalgebra $S(\h)$ by giving a
description of the subspace $\nil(\trigValm)$
of $\trigValm$ in terms of the basis $\{\vec_T\}_{T\in\TabRC(\alm)}$.
\def\lmcl{{\alm[0]}}

Define $\polyH$ as the subalgebra of $\trigH$ generated by
$\F\Preg$, $\F\W$ and $S(\h)$:
$$\polyH=\F\Preg\cdot\F\W \cdot S(\h),$$
which is also a subalgebra of $\ratH$.
% and $\trigH$.
%Note that it is also a subalgebra of $\ratH$. %$\polyH\subset \ratH$.
%Let $\affWreg$ denote
%the semigroup $\Preg\rtimes W\subset \affW$ as before.

We treat a general skew diagram
 $\diag\in\dset_\peri^n$ for a while.
Identify the $\affH$-module $\affV(\diag)$
with the subspace
$\bigoplus_{T\in 
\tabRC({\diag})}\F {v}_T$ of $\trigVdiag$.
Let $\trigVdiag_{+}$ be the $\polyH$-submodule of
$\trigVdiag$ generated by $\affV(\diag)$:
%an $\polyH$-module
$$\trigVdiag_{+} =\polyH \cdot\affV(\diag).$$ %\subseteq \trigVdiag$.
Define
\begin{align}
%\Tab_\peri(\adiag)_{+}&=
%\{T\in\Tab_\peri({\adiag})\mid T(u)\in\Z_{\geq1}\ 
%\ \forall u\in\diag\},\\
 \TabRC(\adiag)_{+}&=
\{T\in\TabRC(\adiag)\mid T(u)\in\Z_{\geq1}\
 \ \forall u\in\diag\}.
\end{align}
\def\T{S}
%%%%%%%%%%%%%%%%%%%%%%%%%%%%%%
\begin{lemma}\label{lem;polyV=tabp}
As a subspace of $\trigVdiag=\F\TabRC(\adiag)$, it holds that
$$\trigVdiag_{+}=\F \TabRC(\adiag)_{+}.
%\bigoplus_{T\in \TabRC(\adiag)_{+}}\F \vec_T.
$$
\end{lemma}
%%%%%%%%%%%%%%%%%%%%%%%%%%%%%%%%
\noindent{\it Proof.}
Let $T\in\TabRC(\adiag)_{+}$.
If
$s_1 T,\dots,s_{n-1}T$ and $\pi T$ are
in $\TabRC(\adiag)$, then
they are actually in $\TabRC(\adiag)_{+}$. 
By the formulas in Theorem \ref{th;tabrep_trig}-(i), 
it follows that $\vec_{s_1T},\dots,\vec_{s_{n-1}T},$ and $\pi \vec_T$
are in $\bigoplus_{T\in \TabRC(\adiag)_{+}}\F \vec_T
=\F \TabRC(\adiag)_{+}$.
Therefore $\F \TabRC(\adiag)_{+}$
is closed under the action of $\affWreg$.
Since  $\F \TabRC(\adiag)_{+}.$
includes $\affV(\diag)=\bigoplus_{T\in \tabRC({\diag})}\F \vec_T$, we have
$\trigVdiag_{+}=\F \affWreg\cdot \affV(\diag)
\subseteq \F \TabRC(\adiag)_{+}.$

\def\Zlm{{{\WRC}_S}}
\def\Zlmp{{\WRC_S^\circ}}

Let us see the opposite inclusion.
Fix any $\T\in\tabRC(\diag)\subseteq \TabRC(\adiag)$ and put
\begin{align*}
\Zlm&=\{w\in\affW\mid w\T\in\TabRC(\adiag)\},\ \
\Zlmp=\{w\in\affW\mid w\T\in\TabRC(\adiag)_{+}\}.
\end{align*}
Observe that $\Zlm\cdot  S=\TabRC(\adiag)$,
$\Zlmp\cdot S=\TabRC(\adiag)_{+}$ and
$\Zlmp=\{w\in \Zlm
\mid w(i)\geq 1\ \forall i\in[1,n]\}$.  

We will prove that $\vec_{w\T}\in\trigVdiag_{+}$
for every $w\in \Zlmp$ by induction on the length of $w$.

Let $w\in\Zlmp$ with  $l(w)=0$.
Then $w=\pi^p$ for some $p\in\Z$
and it is easy to see that
% $w\in\Zlmp$ implies
$p\in\Z_{\geq0}$.
Assume that the statement is true for all $w\in\Zlmp$ with
$l(w)<k$.

Let $w\in\Zlmp$ with $l(w)=k$.
Take $i\in[0,n-1]$ and $x\in\affW$ such that $w=s_ix$ and
$l(w)=l(x)+1$.
Then  $x\in \Zlm$ by Lemma~\ref{lem;weakorder}.

Suppose that $i\in[1,n-1]$.
Then it follows that $x=s_iw\in\Zlmp$
and hence $\vec_{x\T}\in \trigVdiag_{+}$ by
 induction hypothesis.
We have
$$\vec_{w\T}=\vec_{s_ix\T}
=(1+\aa)^{-1}(\vec_{x\T}+\aa s_i\vec_{x\T})
\in \F \affWreg \vec_{x\T},$$
where 
 $\aa:=\con_{x\T}(i)-\con_{x\T}(i+1)$,
 and $1+\aa\neq 0$ by Lemma~\ref{lem;noncritical}.
Therefore $\vec_{w\T}\in  \trigVdiag_{+}$.

Suppose that $i=0$.
As a step, we will show that $\pi^{-1}x\in\Zlmp$.
If $\pi^{-1}x\notin\Zlmp$ then there exists $j\in[1,n]$
such that
$\pi^{-1}x(j)\in\Z_{\leq0}$.
This means $x(j)=1$  since $x(j)\in \Z_{\geq1}$ 
and $\pi^{-1} x(j)=x(j)-1$.
Therefore $w(j)=s_0x(j)=0\notin[1,n]$.
This contradicts $w\in \Zlmp$.
Therefore we have $\pi^{-1}x\in \Zlmp$. 
% and hence $\pi^{-1}x\in\affWreg$.
%We have
%$w=s_0 x=\pi s_{n-1} (\pi^{-1} x)\in\$.

Put $y=\pi^{-1}x$.
The induction
hypothesis implies $\vec_{y\T}\in\trigVdiag_{+}$.

%It easily follows that
Since $y\in\Zlm$ and  $s_{n-1}y=\pi^{-1}w\in \Zlm$,
we have $\aa:=\con_{y\T}(n-1)-\con_{y\T}(n)\neq -1$.
Therefore we have
$$\vec_{w\T}=\vec_{\pi s_{n-1} y\T}=
(1+\aa)^{-1}(\pi \vec_{y\T}+\aa \pi s_{n-1}\vec_{y\T})
\in \F \affWreg \vec_{y\T},$$
and hence $\vec_{w\T}\in \trigVdiag_{+}$.
Therefore we have
$\trigVdiag_{+}
\supseteq\bigoplus_{T\in \TabRC(\adiag)_{+}}\F \vec_T=\F \TabRC(\adiag)_{+}$,
and the statement is proved.
\qed
%
%%%%%%%%%%%%%%%%%%%%%%%%%%%%%%%
\begin{proposition}
\label{pr;positivesubspace}
Let $\kappa\in\Qgeq$ and let $\lm\in
\part_\kappa(\m,n)$. Then
$$ \nil(\trigValm)
=\F \TabRC(\adiag)_{+}$$
%\bigoplus_{T\in\TabRCplm}\F \vec_T,$$
as a subspace of $\trigValm=\F \TabRC(\adiag)$.
\end{proposition}
\noindent{\it Proof.}
It is enough to prove $ \nil(\trigValm)=\polyVlm$.
Since $\affV(\lm)\subset\nil(\trigValm)$,
it holds that $\polyVlm\subseteq  \nil(\trigValm)$.
Since $ \nil(\trigValm)$ is an irreducible $\ratH$-module,
it is enough to show that $\polyVlm$ is closed under
the action of $\ratH$.
Put
$\omega=\sum_{i\in[1,n]}\xch_i\in\ratH\subset\trigH$.
Then
\begin{align*}
 [\omega,\ech_i]&=\kappa\xch_i,\quad [\omega,x_i]=\kappa\quad %(i\in[1,n]).
%\label{eq;omega_e}
(i\in[1,n]),\\ %\label{eq;omega_ech}\\
[\omega, s_i]&=0\quad (i\in [1,n-1]). %\label{eq;omega_s}
\end{align*}
This implies that $\omega$ and the elements in $\polyH$
generate the algebra $\ratH$.
It is easy to check that  $\omega v=0$ for all $v\in\affV(\lm)$.
Hence it follows
% from \eqref{eq;omega_s} and \eqref{eq;omega_e}
that  $\omega$ preserves $\polyVlm =\F\affWreg\cdot \affV(\lm)
=\F\Preg\*\F\W\cdot\affV(\lm)$.
Therefore $\ratH \cdot\polyVlm\subseteq\polyVlm$.
%The latter part of the statement follows dorectly from
%Theorem~\ref{th;tabrep_trig}.
\qed

\smallskip
Combined with Proposition \ref{pr;L=nil}, we obtained
\begin{equation}
  \label{eq;tabreal}
  \ratL(\lm)\cong \F \TabRC(\adiag)_{+}.
\end{equation}
%%%%%%%%%%%%%%%%%%%%%%%%%%%%%%%%%%%
\subsection{Consequences}\label{ss;consequences}
%%%%%%%%%%%%%%%%%%%%%%%%%%%%%%%%%%%%%
%
Let $\kappa\in\Qgeq$ and $m\in\Z_{\geq1}$, 
and put $\peri=\mkappa$. 
Let $\lm\in\part_\kappa(m,n)$.
As a consequence of the realization \eqref{eq;tabreal}
of $\ratL(\lm)$,
%Proposition \ref{pr;positivesubspace},
we obtain the following:
\begin{theorem}\label{cor;tabrep_rat}
%Let $\kappa\in\F\setminus\Q_{\leq 0}$ and 
The $\ratH$-module $\ratVlm$ admits the following weight decomposition
with respect to the subalgebra $S(\h)=\F[\underline{\e}^\vee]:$
%Then 
$$\ratVlm=\bigoplus_{T\in\TabRCplm}
\ratVlm_{\weight_T}.$$
Moreover, $\dim_\F \ratVlm_{\weight_T}=1$ for all 
$T\in\TabRCplm$.
\end{theorem}
%%%%%%%%%%%%%
%%%%%%%%%%%%%%%%%%%%%%%%%%%%%%%%%%%%%%%%%%%%%%%%%%%%%%%%%%%%%
Observe that the degenerate affine Hecke algebra $\affH=S(\h)\F\W
\subset\trigH$ is a subalgebra of $\ratH$ under our identification
$\ratH\subset \trigH$.
As a consequence of Theorem \ref{th;restriction_trig}, we have
\begin{theorem}
%Let $\lm\in\part_\kappa(m,n)$. Then
%\smallskip\noindent
%$\mathrm{(i)}$
%
The $\ratH$-module $\ratVlm$ admits the following decomposition
as an $\affH$-module$:$
$$\ratVlm\!\downarrow_{\affH}
=\bigoplus_{\p\in\PPp(\lm)}
\affV(\p^{-1}(0)).$$
%where $\hat\lm^\p=\p^{-1}(0)$.
\end{theorem}
Recall that $\chi_\zeta:S(\h)^\W\to\F$ denotes
the character corresponding to $\zeta\in \h^*$.
We denote $\chi_\p=\chi_{\weight_T}$ for $\p\in\PP_\peri(\diag)$ as before,
where $T$ is a standard tableau on $\adiag$
%n element of $\TabRC(\adiag)$ 
such that $\projmapop{T}=\p$.
By Theorem \ref{th;VandPP}, we have
%%%%%%%%%%%%%%%%%%%%%%%%%%%%%%%%%%%%%%%%%%%%%%%%%%%%
\begin{theorem}\label{th;VandPP_rat}
%Let $\lm\in\part_\kappa(m,n)$. Then
The $\ratH$-module $\ratVlm$ admits the following decomposition
as an $S(\h)^\W$-module$:$
$$\unip_+\ratVlm=\bigoplus_{\p\in\PPCp(\lm)}
\unip_+\ratVlm^{\chi_\p},$$
%\smallskip\noindent
%$\mathrm{(iii)}$
$$\unip_-\ratVlm=\bigoplus_{\p\in\PPRp(\lm)}
\unip_-\ratVlm^{\chi_\p}.$$
Moreover,  $\dim_\F\; \unip_+\ratVlm^{\chi_\p}=1$
for all $\p\in\PPCp(\lm)$, and
 $\dim_\F\; \unip_-\ratVlm^{\chi_\p}=1$
for all $\p\in\PPRp(\lm)$.
\end{theorem}
By Proposition \ref{pr;tabres_PPC} and 
Theorem \ref{th;VandPP_rat},
we obtain the following:
%%%%%%%%%%%%%%%%%%%%%%%%%%%%%%%%%%%%%%%%%%%%%%%
\begin{theorem}\label{th;sp_decomp_rat}
The $\ratH$-module $\ratVlm$ admits the following decomposition
as an $S(\h)^\W$-module$:$
  $$\unip_+\ratVlm
=\bigoplus_{\zeta\in \Preg^-, \ T\in\restab_{\peri}(\lm)}
%{\zeta\in\dst^{\mathrm{res}}_\circ(\alm)}
\unip_+\ratVlm^{\chi_{\zeta+\ccwt_T}},$$
%where 
%$$\dst^{\mathrm{res}}_\circ(\alm)=
%\left\{t_\zeta T_\sharp\in\Tab(\alm)\mid
%\zeta\in P_\circ^-,\ T\in\restab_{\peri}(\lm)\right\}.$$
%$\subseteq \TabRC(\adiag).$
Moreover, 
$\unip_+\ratVlm^{\chi_{\zeta+\ccwt_T}}
=\F\, \unip_+ \vec_{t_{\zeta+\ccwt_T} T}$
for all $\zeta\in \Preg^-$
and $T\in\restab_\pri(\lm)$.
\end{theorem}
%%%%%%%%%%%%%%%%%%%%%%%%%%%%%%%%%%%%%%%%%%%%%%%%%%%%%%
\begin{remark}
%Suppose that $\kappa\notin$
%(i)
 Put $C=\Q_{\geq0}\cup\{-2/(2k+1)\mid k\in\Z_{\geq0}\}$.
%Moreover suppose that 
%$\kappa\notin\{ \frac{2}{2r-1}\mid r\in\Z_{\geq1}\}$.
Then it is shown by Gordon-Stafford \cite{GS1} 
(see also \cite{BEG;spherical}) that
 the assignment $M\mapsto \unip_+M$
gives a Morita equivalence
 $\ratH$-$\mod\to
\unip_+\ratH\unip_+$-$\mod$ if $\kappa\notin C$.
Via the isomorphism  $\ratH\cong\H_{-\kappa}$
in Lemma \ref{lem;invo},
it also holds that
 the correspondence $M\mapsto \unip_-M$
gives a Morita equivalence
 $\ratH$-$\mod\to
\unip_-\ratH\unip_-$-$\mod$
if  $-\kappa\notin C$
%$\Q_{\leq0}\cup\{2/(2k+1)\mid k\in\Z_{\geq0}\}$.

The algebra $\unip_+\ratH\unip_+$ is called the 
{\it spherical subalgebra} of $\ratH$.

It should be mentioned that
parallel statements can be shown for $\trigH$-modules
using the isomorphism $\XL\otimes_{\X}\ratH\cong\trigH$.
\end{remark}
%%%%%%%%%%%%%%%%%%%%%%%%%%%%%%%%%%%%%%%%%%%%%%%%%%%%%%%%%%%%%%
\section{Characters}
%%%%%%%%%%%%%%%%%%%%%%%%%%%%%%%%%%%%%%%%%%%%%%%%%%%%%%%%%%%%%%%
We introduce characters for
 $\unip_\pm\ratVlm$, and compute them
 in several cases using
the results in \S \ref{sec;generating_function}
and \S \ref{ss;consequences}.

%%%%%%%%%%%%%%%%%%%%%%%%%
\subsection{Characters and generating functions}
%%%%%%%%%%%%%%%%%%%%%%%%%

The algebra $\ratH$ has the grading operator
\begin{equation}
\label{eq:der}
\der
=\kappa^{-1}\sum_{i=1}^{n}\ech_i
=\kappa^{-1}\sum_{i=1}^{n}x_i y_i+
{\kappa^{-1}}\sum_{1\leq i<j\leq n}s_{ij}.
\end{equation}
%where
%$\der'=\kappa^{-1}\sum_{i=1}^{n}x_i y_i$ and $\der_{\W}=
%{\kappa^{-1}}\sum_{1\leq i<j\leq n}s_{ij}$,
The element $\der$
 belongs to $S(\h)^\W$ and satisfies 
\begin{alignat*}{2}
&[\der,{x_i}]={x_i}, 
\ [\der,\ech_i]=-\ech_i\quad &(i\in[1,n]),\\
& [\der,w]=0\quad &(w\in\W).
\end{alignat*}
In particular, 
 $\der$ preserves the  $\W$-(anti-)invariant subspace 
of any $\ratH$-module.

We put $\der'=\kappa^{-1}\sum_{i=1}^{n}x_i y_i$ and 
$\der_{\W}={\kappa^{-1}}\sum_{1\leq i<j\leq n}s_{ij}$.
Then $\der=\der'+\der_{\W}$ and $[\der',\der_{\W}]=0$.

For an $\F[\der]$-module $M$ and $d\in\F$, put
$M{[d]}=\{v\in M\mid (\der-d)^k v=0 \hbox{ for }k>>1\}$.
%Let $\lm\in\Lambda^+_\kappa(\m,n)$
%and put
%$\unip\ratVlm_{(k)}=\{v\in \unip\ratVlm\mid\der v=kv\}$ 
%for $k\in\F$.
For  $M\in\O(\ratH)$, it follows that $M$ is finitely generated 
over $\F[\udl x]$,
and hence we have $\dim_\F M{[d]}<\infty$, and the trace
%we define
%and $\dim_\F \unip M<\infty$ 
%It follows that 
$$\Ch M=\sum_{d\in\F}(\dim_\F M[d])q^d$$
is defined.
%where $\mathrm{Tr}(M;F)$ denotes the trace of the operator $F$ on
%$M$.
Similarly, we can define $\Ch{\unip_\pm M}$,
which we call %$\Ch{\unip_+ M}$ %(resp.,  and $\Ch{\unip_- M}$ )
{\it the spherical character} 
%(resp., {\it the co-spherical character}) 
of $M$.

Let $\lm\in \Lambda^+_\kappa(m,n)$. Put
%\begin{align*}
%&|T|=\sum_{u\in\lm}T(u)\ \hbox{ for }T\in\Tab(\alm),\\&
$|\p|=\sum_{u\in\lm}\p(u)$ for $\p\in\PP(\lm)$
as before.
%\end{align*}
The following can be checked easily.
%%%%%%%%%%%%%%%%%%%%%%%%%%%%%%%%%%%%%%%%%%%%%%%%
\begin{lemma}
Let $T\in\TabRC(\alm)$.
Then for $\vec_T\in\trigValm$,
\begin{equation}
\der \vec_T=(|\projmapop{T}|+\cas_\lm)\vec_T,
\end{equation}
where 
\begin{equation}\label{eq;cas}
\cas_\lm=\frac{1}{2\kappa}\sum_{i=1}^\m
\lm_i(\lm_i-2i+1).
\end{equation}
\end{lemma}
%%%%%%%%%%%%%%%%%%%%%%%%%%%%%%%%%%%%%%%%%%%%%%%%%%%%%
By Theorem \ref{th;VandPP_rat},
the computation of the characters are reduced to the computation 
of the generating functions for plane partitions:
\begin{corollary}\label{cor;char_GF}
Let $\kappa\in\Qgeq$ and $m\in\Z_{\geq1}$, and put $\peri=\mkappa$. 
Let $\lm\in\Lambda^+_\kappa(m,n)$. Then

\begin{align*}
%$\mathrm{(i)}$
%
%&\Ch{\ratVlm} 
%=q^{\cas_\lm}\sum_{\p\in\trigPPp(\lm)_{\geq0}}
%q^{|\p|},\\
&\Ch{\unip_+\ratVlm}=
q^{\cas_\lm}\sum_{\p\in\PPCp(\lm)}
q^{|\p|},\\
&\Ch{\unip_-\ratVlm}=
q^{\cas_\lm}\sum_{\p\in\PPRp(\lm)}
q^{|\p|}.
\end{align*}
\end{corollary}
%
%$\kappa$ is ``generic''.

%We also compute their spherical characters when
%\smallskip\noindent(ii)
% $\kappa$ is an integer.

%%%%%%%%%%%%%%%%%%%%%%%%%%%%%%%%%%%%%%%%%%%%%%%%%%%%%%%%%%%%%%
\subsection{Spherical characters of irreducible modules
and Kostka polynomials}
%%%%%%%%%%%%%%%%%%%%%%%%%%%%%%%%%%%%%%%%%%%%%%%%%%%%%%%%%%%%%%
Combining Corollary \ref{cor;char_GF} and the results
in Section \ref{sec;generating_function},
we obtain character formulas in the following cases:

\smallskip\noindent
(i) $\kappa\in\Z_{\geq1}$.

\smallskip\noindent
(ii) $\kappa$ is "generic".

\medskip
%%%%%%%%%%%%%%%%%%%
First, let $\kappa\in\Z_{\geq1}$.
By Theorem \ref{th;sp_decomp_rat} and Theorem \ref{th;GFCperi},
 we obtain the following
formula, which is
 announced in the previous paper  by the author \cite{Su;conformal}.
%%%%%%%%%%%%%%%%%
\begin{theorem}\label{th;spchar_rat}
Let $\m,\kappa \in\Z_{\geq1}$
and
 $\lm\in \Lambda^+_{\kappa}(\m,n)$ .
Then
  \begin{equation*}
    \Ch{\unip_+\ratVlm}=
\frac{q^{\cas_\lm}}{\qfac{n}}
%{\mathop{\Pi}\limits_{i=1}^n(1-q^i)}
{\check K}_{\lm\; (1^n)}^{(\kappa-\m)}(q),
%=
%\frac{q^{\cas_\lm+\frac{1}{2}n(n-1)}}{\qfac{n}}
%{\mathop{\Pi}\limits_{i=1}^n(1-q^i)}
%{K}_{\lm\; (1^n)}^{(\kappa-\m)}(q^{-1}).
  \end{equation*}
where 
  $\qfac{n}=(1-q)(1-q^2)\dots (1-q^n)$.
\end{theorem}
\begin{remark}
 The formula for $\unip_-\ratVlm$ with $\kappa\in\Z_{\leq-1}$
 can be obtained through the isomorphism 
$\ratH\cong {\mathcal H}_{-\kappa}$.
 But, the formula for $\unip_-\ratVlm$ with $\kappa\in\Z_{\geq1}$
has not been obtained (see also Remark \ref{rem;row}).
%It might not be as simple as the one for  $\unip_+\ratVlm$.
It seems remarkable that $\unip_+\ratVlm$
has such a simple character formula.
\end{remark}
%%%%%%%%%%%%%%%%%%%%%

Next, we consider the case where $\kappa$ is generic.
Let $\kappa\in\F\setminus\{0\},\ m\in\Z_{\geq1}$, and let
$\diag\in\dset_\peri^n$ with $\peri=\mkappa$.
 \begin{definition}\label{def;generic}
The $\kappa$ is said to be {\it generic} for $\diag$
if %the following condition holds:
%\smallskip\noindent
$$\{(a+k+1,b+k)\in\F\times\Z 
\mid k\in\Z_{\geq0}\}\cap\adiag_\peri
\subseteq\diag
\hbox{ for any }(a,b)\in\diag,$$
where $\peri={(-m,\kappa-m)}$.
 \end{definition}
It is obvious that for $\diag\subset\Z\times\Z$, 
irrational $\kappa$ is generic. 
%More precisely, we have
The following can be shown easily.
%%%%%%%%%%%%%%%%%%%%%%%%%%%%%%%%%%%%%%%%%%%%%%%%%%%%%%%%%%%%%%%%%%%%%%%%%%%%%
\begin{proposition}\label{pr;generic}
Let $\kappa\in\F\setminus\Q_{\leq0}$ and
$\lm\in\part_\kappa(m,n)$.
Then $\kappa$ is generic for $\lm$ if and only if one of the following holds$:$

\smallskip\noindent
$(a)$ $\kappa\notin\Q$.

\smallskip\noindent
$(b)$ $\kappa=s/r$ with $r,s\in\Z_{\geq1}$, $(s,r)=1$, and $\lm_1\geq s-m+1$.
\end{proposition}
%%%%%%%%%%%%%%%%%%%%%%%%%%%%%%%%%%%%

The following is obvious
\begin{lemma}
  If $\kappa$ is generic for $\diag$, then
$\PP(\diag)=\PPp(\diag)$, $\PPR(\diag)=\PPRp(\diag)$
and $\PPC(\diag)=\PPCp(\diag)$.
\end{lemma}
%%%%%%%%%%%%%%%%%%%%%%%%%%%%%%%%%%%%%%%%%%%%%%%%%%%%%%%%%%%
As a direct consequence of Proposition 
\ref{pr;char_classical}, we have 
the following:
%%%%%%%%%%%%%%%%%%%%%%%%%%%%%%%%%%%%%%%%%%%
\begin{corollary}\label{cor;char_generic}
Let $\kappa\in\F\setminus\{0\}$, $\m\in\Z_{\geq1}$ and
$\lm\in\Lambda^+(\m,n)$. 
Suppose that $\kappa$ is generic for $\lm$.
Then $$\Ch {\unip_+ \ratVlm}
=\frac{q^{\cas_\lm}}{\qfac{n}}
{\check K}_{\lm\; (1^n)}(q),$$
$$\Ch {\unip_- \ratVlm}
=\frac{q^{\cas_\lm}}{\qfac{n}}
{K}_{\lm\; (1^n)}(q).$$
\end{corollary}
%%%%%%%%%%%%%%%%%%%%%%%%%%%%%%%%%%%%%%%%%%%%%%%%%%%%%%%%%%%%%%%%

%%%%%%%%%%%%%%%%%%%%%%%%%%%%%%%%%%%%%%%%%%
\subsection{Characters of standard modules and their irreducibility}
%{Character of standard modules}
%%%%%%%%%%%%%%%%%%%%%%%%%%%%%%%%%%%%%%%
Let $\lm\vdash n$.
We have an isomorphism
 $$\ratst(\lm)\isomto \X\otimes S_\lm$$
 as an $\F\W$-module.
Observe that it is also an isomorphism of $\F[\der']$-modules,
where $\der'$ acts on $\X\otimes S_\lm$ %the right hand side as
as  $\sum_{i\in[1,n]}x_i\frac{\der}{\der x_i}\otimes \id_{S_\lm}$.
The following is straightforward:
%%%%%%%%%%%%%%%%%%%%%%%%%%%%%%%%%%%%%
\begin{proposition} 
Let $\kappa\in\F\setminus\{0\}$ and
let $\lm\vdash n$. Then
$$\Ch{\ratst(\lm)}=%\Ch{\ratVlm}=
\frac{q^{d_\lm}}{(1-q)^n}\dim_\F S_\lm.$$
\end{proposition}
%%%%%%%%%%%%%%%%%%%%%%%%%%%%%%%%%%%%%5
%Using Corollary \ref{cor;char_generic},
We will see that the genericity of $\kappa$ is
a sufficient condition for a standard module to be irreducible.
%and will obtain their character formulas.
%

\def\Tmap{\Sigma}
%{{\mathcal{S}}}
For $\zeta=\sum_{i\in[1,n]}\zeta_i\e_i\in P$ and $T\in\tabRC(\diag)$,
define  $j_1<j_2<\dots<j_n$  to be the increasing sequence in $\Z$
obtained as a rearrangement of
$1+n\zeta_1,2+n\zeta_2,\dots,n+n\zeta_n$,
and define the  tableau $\Tmap=\Tmap_{(\zeta,T)}$ on $\adiag$ by
$$\Tmap(u+k\peri)=j_{T(u)}-kn\quad
(u\in\diag,\ k\in\Z).$$
Then it follows from the definition
that $\Tmap$ is row-column increasing on $\diag$.
Suppose now that $\kappa$ is generic for $\diag$. 
Then it is easy to see that
$\Tmap\in\TabRC(\adiag)$,  and hence
 the assignment $(\zeta,T)\mapsto\Tmap_{(\zeta,T)}$ gives a map
$P\times\tabRC(\diag)\to\TabRC(\adiag)$.
%%%%%%%%%%%%%%%%%%%%%%%%%%%%%%%%%%%%%%%%%%%%%%%%%%%%%%%%%%%%%
\begin{lemma}\label{lem;Ptab_Tab}
Let $\diag\in\dset^n_\peri$ 
and suppose that $\kappa$ is generic for $\diag$.
Then the correspondence $(\zeta,T)\mapsto \Tmap_{(\zeta,T)}$
above gives bijections
$$P\times\tabRC(\diag)\isomto\TabRC(\adiag),
\quad \Preg\times\tabRC(\diag)\isomto\TabRC(\adiag)_{\posi}.
$$
\end{lemma}
%%%%%%%%%%%%%%%%%%%%%%%%%%%%%%%%%%%%%%%%%%%%%%%%%%%%%%%%%%%%%%
\noindent{\it Proof.}
For $S\in\TabRC(\adiag)$, define
$l_1<\dots<l_n$ as the increasing sequence in $\Z$
obtained as a rearrangement of
$\{S(u)\mid u\in\diag\}$.
Write $l_i=\bar{l}_i+k_in$ with $\bar{l}_i\in[1,n]$ and $k_i\in\Z$,
and define $\zeta_S=\sum_{i\in[1,n]}k_i\e_{\bar{l}_i}\in P$.
Define $T_S\in\tab(\diag)$ by
$$T_S(u)=i\hbox{ for }u\in\diag\hbox{ such that }S(u)=l_i.$$
Then $T_S\in\tabRC(\diag)$ as $S\in\TabRC(\adiag)$.
It is easy to check that the map $\TabRC(\adiag)\to P\times\tabRC(\diag)$
given by the correspondence
$S\mapsto (\zeta_S,T_S)$ is the inverse map of $\Tmap$.
\qed
%%%

%\medskip

%
%%%%%%%%%%%%%%%%%%%%%%%%%%%%%%%%%%%%%%%%%%%%
\begin{proposition}\label{pr;generic_irred}
 Let $\lm\vdash n$.
%$\in\part_\kappa(m,n)$.
 If $\kappa$ is generic for $\lm$.
Then $\ratst(\lm)$ is irreducible and isomorphic to $\ratVlm$.
\end{proposition}
\noindent{\it Proof.}
Observe that for $\Tmap=\Tmap_{(\zeta,T)}$
in Lemma \ref{lem;Ptab_Tab}, we have
$\der \vec_\Tmap=(|\zeta|+\cas_\lm) \vec_\Tmap$,
where $|\zeta|=\sum_{i\in[1,n]}\zeta_i$.
Hence the lemma implies
$$\Ch{\ratVlm}=\dim_\F S_\lm\cdot
q^{\cas_\lm}\sum_{\zeta\in \Preg}q^{|\zeta|}
=\frac{q^{\cas_\lm}}{(1-q)^n}\dim_\F S_\lm .$$
Because we have a surjective homomorphism 
$\ratst(\lm)\to\ratVlm$
and both sides have the same character,
they are isomorphic.
\qed

%\medskip

%%%%%%%%%%%%%%%%%%%%%%%%%%%%%%%%%%%%%%%%%%%%%%%%%%%%%%%%%%%%
\begin{proposition}
$($cf. \cite[Proposition 3.9]{GS1}$)$
\label{pr;sphericalchar_st}
  Let $\kappa\in\F\setminus\{0\}$, $\m\in\Z_{\geq1}$ and
$\lm\in\Lambda^+(\m,n)$. Then
$$\Ch {\unip_+ \ratst(\lm)}
=\frac{q^{\cas_\lm}}{\qfac{n}}
%{\mathop{\Pi}\limits_{i\in[1,n]}(1-q^i)}
{\check K}_{\lm\; (1^n)}(q),$$
$$\Ch {\unip_- \ratst(\lm)}
=\frac{q^{\cas_\lm}}{\qfac{n}}
%{\mathop{\Pi}\limits_{i\in[1,n]}(1-q^i)}
{K}_{\lm\; (1^n)}(q).$$
%where 
%$\cas_\lm=\frac{1}{2\kappa}\sum_{i=1}^\m
%\lm_i(\lm_i-2i+1).$
\end{proposition}
%%%%%%%%%%%%%%%%%%%%%%%%%%%%%%%%%%%%%%%%%%%%%%%%%%%%%%%%%%%%%%%%
%
\noindent{\it Proof.}
Note that $\Ch{\unip_\pm\ratst(\lm)}=q^{\cas_\lm}
\mathrm{Tr}(\unip_\pm\ratst(\lm);q^{\der'})$
and $\mathrm{Tr}(\unip_\pm\ratst(\lm);q^{\der'})$ does not depend on 
$\kappa$. 
Hence it is enough to show the formula  
when $\kappa$ is generic for $\lm$ (for example, $\kappa\notin\Q$).
%, and
%compute $\Ch{\unip_+\ratVlm}=\Ch{\unip_+\ratst(\lm)}$.
In this case, we have
 $\Ch{\unip_\pm\ratVlm}=\Ch{\unip_\pm\ratst(\lm)}$,
 and the statement follows from Corollary \ref{cor;char_generic}.
\qed

\begin{remark}
Proposition \ref{pr;sphericalchar_st} can be also derived from 
the result by Garsia-Procesi \cite{GP}.
As is well-known,
 $\X$ is a free module over $\X^\W$,
and furthermore
$\X\cong\X^\W\otimes \mathcal{R}$ as a graded $\F\W$-module, 
where $\mathcal{R}$ (called the coinvariant algebra)
is by definition the quotient algebra of $\X$ by
the ideal generated by the elementary symmetric functions
of positive degree.

In \cite{GP}, it is proved that in the decomposition 
\begin{equation}
  \label{eq;coinvalg}
\mathcal{R}=\oplus_{\lm\vdash n}M_\lm\otimes S_\lm
\end{equation}
as a graded $\F\W$-module, 
the graded dimension of $M_\lm$ %$\mathrm{Tr}(M_\lm;q^{\der'})$
 is given by %the (deformed) Kostka polynomial
$\check{K}_{\lm\, (1^n)}(q)$.
This implies the formula in Proposition \ref{pr;sphericalchar_st}.

%It should be also remarked that 
Conversely, our proof of
Proposition \ref{pr;sphericalchar_st}
gives an alternative derivation of the formula 
$\mathrm{Tr}(M_\lm;q^{\der'})=\check{K}_{\lm\, (1^n)}(q)$
in \cite{GP}.
%Since we have $\unip_+ \ratst(\lm)=(\X\otimes S_\lm)^W
%=\X^\W\otimes M_\lm$, 
%the statement follows from $\der=\der'+\cas_\lm$.
\end{remark}
%%%%%%%%%%%%%%%%%%%%%%%%%%%%%%%%%%%%%%%%%%%%%
\subsection{Characters for  single column representations}
\label{ss;rep_single_column}
%%%%%%%%%%%%%%%%%%%%%%%%%%%%%%%%%%%%%%%%%%%%%%%
%
%%%%%%%%%%%%%%%%%%%%%%%%%%%%%%%%%%%%%%%%%%%%%%%%%%%%%
\def\Pset{M}
%In this section, we
As a result for 
non-generic and non-integral $\kappa$, 
we give character formulas
for  the tableaux representation corresponding to 
the single column
$\vert=\{(a,1)\in\Z\times\Z \mid a\in[1,n]\}$
using Proposition \ref{pr;GF_single_column}.
%Observe that $\vert\in\dset_\peri^n$ with $\peri=(n,n-n/r)$.
%e consider the structure of the $\trigH$-module $\trigV(\avert)$.
%Let $\kappa\in\Qgeq$.
\begin{lemma}\label{lem;unless_n/r}
Let $\kappa\in\F\setminus\Q_{\leq0}$ and put $\peri=(-n,\kappa-n)$.

\smallskip\noindent
$\mathrm{(i)}$ $\vert\in\dset_\peri^n$
%$\part_\kappa(n,n)$
if and only if $\kappa\notin \Q$
or $\kappa=s/r$ with
$s\in\Z_{\geq n}$, $r\in\Z_{\geq1}$ and $(s,r)=1$.

\smallskip\noindent
$\mathrm{(ii)}$ Suppose $\vert\in\dset_\peri^n$.
Then $\ratst(\vert)\cong\ratV(\vert)$ unless
$\kappa=n/r$ with
$r\in\Z_{\geq1}$ and $(n,r)=1$.
\end{lemma}
\noindent{\it Proof.}
Note that $\vert\in\dset_\peri^n\Leftrightarrow \vert\in\part_\kappa(n,n)$,
and the statement (i) is obvious.

By Proposition \ref{pr;generic}, $\kappa$ is generic 
%(in the sense of Definition \ref{def;generic}) 
unless $\kappa=n/r$.
Hence the statement (ii) follows from 
 Proposition \ref{pr;generic_irred}.
%,which we will prove later. 
\qed

\medskip
%%%%%%%%%%%%%%%%%%%%%%%%%%%%%%%%%%%%%%%%%%%%%%%%%%%%%%%%%%%%%%%%
In the rest, we suppose 
$\kappa=n/r$ with $r\in\Z_{\geq1}$, $(n,r)=1$.
%%%%%%%%%%
\begin{proposition}
%\smallskip\noindent
%  $\mathrm{(i)}$
The $\ratH$-module $\ratV(\vert)$ is 
a free module over $\F[\pi]$ of rank $r^{n-1}$. Moreover,
\begin{align*}
\Ch{\ratV(\vert)}
&=\frac{q^{-\frac{1}{2}r(n-1)}}{1-q}
\left(\frac{1-q^r}{1-q}\right)^{n-1},\\
\Ch{\unip_-\ratV(\vert)}&=
%\frac{q^{-\frac{1}{2}r(n-1)}}{1-q^n}
q^{-\frac{1}{2}r(n-1)}
\frac{ \qfac{n+r-1} }{ \qfac{n} \qfac{r} },\\
\Ch{\unip_+\ratV(\vert)}&=%\frac{q^{-\frac{1}{2}r(n-1)}}{1-q^n}
\begin{cases}
0\ \ &\text{if}\ r<n,\\
q^{-\frac{1}{2}r(n-1)}
\frac{ \qfac{r-1}}{ \qfac{n} \qfac{r-n} }\ \ 
&\text{if }r\geq n.
\end{cases}
\end{align*}
\end{proposition}
%%%%%%%%%%%%%%%%%%%%%%%%%%%%%%%%%%%%%%%%%%%%%%%%%%%%%%%%%%%%%%%%%
\begin{remark}
For the rational Cherednik algebra $\ratH(SL_n)$
of type $SL_n$, formulas for $\Ch{\ratV(\vert)}$ and $\Ch{\unip_-\ratV(\vert)}$
have been obtained by Berest-Etingof-Ginzburg \cite{BEG}
by a different method.
(They furthermore computed 
$\mathrm{Tr}({\ratV(\vert)};wq^\der)$ for any $w\in\W$.)
In the same paper
%\cite{BEG} 
 it is also proved that
finite-dimensional representations for $\ratH(SL_n)$
exist only for $\kappa=\pm n/r$ with $r\in\Z_{\geq1}$, $(r,n)=1$,
and moreover
any finite-dimensional irreducible module 
is isomorphic to  $\ratL(\vert)$ in ``$+$'' case
and is isomorphic to $\ratL(\vert')$ in ``$-$'' case
(see also \cite[p.65]{Ch;fourier}).
\end{remark}
%\begin{remark}
%\end{remark}
%Put 
%\begin{proposition}
%
%The $\trigV(\avert)$ is a free module over $\F[\pi,\pi^{-1}]$
%of rank $r^{n-1}$.
%\end{proposition}
\noindent{\it Proof.}
The formulas for $\unip_\pm\ratV(\vert)$
 are direct consequences of
Proposition \ref{pr;GF_single_column}.
We shall prove the statement for $\ratV(\vert)$.

By similar argument as in the proof of Proposition \ref{pr;GF_single_column},
%Lemma \ref{lem;equi_rstep}, it can be checked that 
it follows that 
\begin{eqnarray*}
%\label{eq;tabrcvert}
\TabRC(\avert)_{+}=
\label{eq;tabrcver}
\left\{T\in\Tab_\peri(\avert)\mid 
%\hbox{ and }
%T\hbox{ satisfies }\eqref{eq;cinc}
0<T(1,1)<T(2,1)<\dots <T(n,1)<T(1,1)+rn
\right\}. 
\end{eqnarray*}

Observe that $\pi\in\affW$ acts on the space 
$\ratL(\vert)=\F\TabRC(\avert)_{+}$ by
 $\pi \vec_T=\vec_{\pi T}$.

Put 
$\TabRC(\avert)_1=\{T\in\TabRC(\avert)\mid T(1,1)=1\}$.
Then it easily follows that 
$\F\TabRC(\avert)_{+}$ is a free module over $\F[\pi]$
with the basis $\{\vec_T \}_{ T\in\TabRC(\avert)_1}$.

Let
$\Pset(n,r)$ denote the set of all maps
from $[2,n]$ to $[0,r-1]$.
For $f\in\Pset(n,r)$, there exists a unique 
$T\in \TabRC(\avert)_1$
 such that %$T(1)=1$
$\{T(2,1).\dots,T(n,1)\}=\{i+f(i)n\mid i\in[1,n]\}$.
From the expression of $\TabRC(\avert)_{+}$, it follows
%\eqref{eq;tabrcvert} 
that this correspondence gives 
a bijection $\Pset(n,r)\to \TabRC(\avert)_1$.
%Put 
%Then we have
%gives a set of the generators
%of $\trigV(\vert)$ over $A$.
In particular, we have 
$\mathrm{rank}_{\F[\pi]} \ratV(\vert)=\sharp \TabRC(\avert)_1=r^{n-1}$.
%\qed

Let $T_f$ denote the element in $\TabRC(\avert)_1$
corresponding to $f \in \Pset(n,r)$ through the bijection
 $\Pset(n,r)\isomto \TabRC(\avert)_1$ above.
Then, we have $\der \vec_{T_f}=(|f|+\cas_{\vert}) \vec_{T_f}$, where
$|f|=\sum_{i\in[2,n]}f(i)$ and  $\cas_{\vert}=-\frac{1}{2}r(n-1)$.
It is easy to see that
$$\sum_{f\in\Pset(n,r)}q^{|f|}=
\left(\frac{1-q^r}{1-q}\right)^{n-1},$$
from which the statement follows.
\qed
%
%%%%%%%%%%%%%%%%%
\begin{appendix}
%%%%%%%%%%%%%%%%%%%%%%%%%%%%%%%%%%%%%%%%%
\section{Classification of irreducible modules 
with weight decomposition}
%%%%%%%%%%%%%%%%%%%%%%%%%%%%%%%%%%%%%%%%%%%%%%%%%
We give a proof of Theorem \ref{th;classification_rat}.
For $\zeta\in\wt$, we define
a function  $F_\zeta$ on $\Z$
by setting 
\begin{equation*}
%\label{eq;extendwt}
F_\zeta(i)=\brac{\zeta}{\ech_{\bar i}}-k\kappa\ \ 
\hbox{ for }i=\bar i+kn \hbox{ with }
\bar i\in[1,n],\ k\in\Z.
\end{equation*} 
To prove Theorem \ref{th;classification_rat}, we use the following lemma:
%%%%%%%%%%%%%%%%%%%%%%%%%%%%%%%%%%%%%%%%%%%%%%%%%%%%%%%%%%%%%%%%
\begin{lemma}$($cf. \cite[Lemma 4.19]{SV}$)$
\label{lem;klemma}
Let $\kappa\in\Qgeq$.
Let $L$ be an irreducible $\trigH$-module 
which belongs to $\Oss(\trigH)$, and let
 $\zeta\in \h^*$ be a weight of $L$.
For any $i,j\in\Z$ such that 
$$ i<j\quad \text{and}\quad \brac{\zeta}{\alch_{ij}}=0,$$
there exist $k_+\in [i+1,j-1]$ and  
$k_-\in[i+1,j-1]$ such that
$$F_\zeta(k_\pm)=F_\zeta(i)\pm1.$$
\end{lemma}
%%%%%%%%%%%%%%%%%%%%%%%%%%%%%%%%%%%%%%%%%%%%%%%%
\noindent{\it Proof.}
The statement has been
proved in  \cite[Lemma 4.19]{SV}
with the restriction $\zeta\in P$ and $\kappa\in\Z$.
The proof can be generalized to our case with
little modification.
\qed
%%%%%%%%%%%%%%%%%%%%%%%%%%%%%%%%%%%%%%%%%%%%%%%%%%%%%%%%%%%%%%%%5

\medskip
\noindent
{\it Proof of Theorem \ref{th;classification_rat}.}

\smallskip
We have seen that  $\ratL(\lm)\cong \nil(\trigValm)$
when $\lm\in\bigsqcup_{\m\in[1,n]}\part_\kappa(\m,n)$, 
and hence $\ratL(\lm)$ is in $\Oss(\ratH)$ 

We suppose that $\lm\in \part(\m,n)\setminus \part_\kappa(\m,n)$, 
and will prove that  $\ratL(\lm)$ does not belong to $\Oss(\ratH)$.
The statement is easily checked when $\kappa\notin\Q$.

Assume that $\kappa\in\Q_{\geq0}$ and write $\kappa=s/r$
with $r,s\in\Z_{\geq1}$, $(r,s)=1$.

\def\rblm{{\mathfrak{t}_\lm}}
Since $S_\lm=\affV(\lm)\subset \ratL(\lm)$,
the content $\con_{\rblm}$
of the row reading tableau $\rblm$ on $\lm$ gives a weight $\weight_{\rblm}$
of $\ratL(\lm)$. 
By the assumption, we have
\begin{equation*}
  \label{ineq1}
s-m-\lm_1+\lm_m\in\Z_{<0}.
\end{equation*}
First, assume that $s< m$.
Put 
%$j=n-\lm_m+m-r$.
$a=m-s$. 
Then $a\in[1,m-1]$ and $(a,1)\in\lm$.
We put
$i=\rblm(m,1)$ and
$j=\rblm(a,1)+rn$.
We have 
$$F_{\weight_{\rblm}}(j)=F_{\weight_{\rblm}}(\rblm(a,1))-r\kappa
=1-a-s=1-m=F_{\weight_{\rblm}}(i).$$

It is easy to see that
\begin{align*}
 F_{\weight_{\rblm}}(k)>F_{\weight_{\rblm}}(i)&\ 
\hbox{ for all }k\in [i+1,n],\\
F_{\weight_{\rblm}}(k)>F_{\weight_{\rblm}}(j)&\ 
\hbox{ for all }k\in [1+rn,j-1],\\
F_{\weight_{\rblm}}(k)\notin\Z&\ \hbox{ for all }k\in [n+1,rn].
\end{align*}
Therefore there are no $k\in[i+1,j-1]$
such that $F_{\weight_{\rblm}}(k)=F_{\weight_{\rblm}}(i)-1$.
By Lemma~\ref{lem;klemma}, $\ratL(\lm)$ is not in $\Oss(\ratH)$.

Next, assume that $s\geq m$.
Put $b=1+\lm_m+s-m.$
Then we have
 $b\in[1,\lm_1]$ and hence $(1,b)\in\lm$.
Put $i={\rblm}(m,\lm_m)=n$
and $j={\rblm}(1,b)+rn$.
Then
$F_{\weight_{\rblm}}(j)=F_{\weight_{\rblm}}(\rblm(1,b))-r\kappa
=(1+\lm_m+s-m)-1 -s=\lm_m-m=
F_{\weight_{\rblm}}(i).$
Similarly to the case $s< m$, it is shown that
there are no $k\in[i+1,j-1]$
such that $F_{\weight_{\rblm}}(k)=F_{\weight_{\rblm}}(i)+1$.
By Lemma~\ref{lem;klemma},  $\ratL(\lm)$ is not in $\Oss(\ratH)$.
\qed
\end{appendix}

\bigskip
\noindent
{\bf Acknowledgment.}
The author would like to thank S. Ariki, S. Okada. and
A. Kirillov
for valuable suggestion and comments.
The author is also grateful to T. Kuwabara and 
M. Kasatani for fruitful discussion.
%%%%%%%%%%%%%%%%%%%%%%%%%%

%%%%%%%%%%%%%%%%%%%%%%%%%%%%%%%%%%%%%%%%%%%%%%%%
\end{document}